# SAMPLE PATH LARGE DEVIATIONS FOR MULTICLASS FEEDFORWARD QUEUEING NETWORKS IN CRITICAL LOADING

By Kurt Majewski

*Siemens AG, Corporate Technology*

We consider multiclass feedforward queueing networks with first in first out and priority service disciplines at the nodes, and class dependent deterministic routing between nodes. The random behavior of the network is constructed from cumulative arrival and service time processes which are assumed to satisfy an appropriate sample path large deviation principle. We establish logarithmic asymptotics of large deviations for waiting time, idle time, queue length, departure and sojourn-time processes in critical loading. This transfers similar results from Puhalskii about single class queueing networks with feedback to multiclass feedforward queueing networks, and complements diffusion approximation results from Peterson. An example with renewal inter arrival and service time processes yields the rate function of a reflected Brownian motion. The model directly captures stationary situations.

**1. Introduction.** A frequently considered asymptotic for analyzing the complex behavior of queueing networks deals with the situation when the queues approach a critical load: Based on a functional central limit theorem for the driving processes, the weak convergence of suitably normalized network behavior to reflected (fractional) Brownian motion can be proved [13, 17, 21, 29, 34]. Here we speak of a *reflected* process when it is restricted to the positive orthant through a Skorokhod map [9, 14]; see Remark 5 after Theorem 2 for a definition. This heavy traffic convergence has been one motivation for a large number of publications about large deviations of reflected Gaussian processes [1, 20, 23, 25, 27, 28]. However, the weak convergence of queueing processes to a reflected Gaussian process, is not appropriate for justifying a relation between large deviations of the queueing processes and large deviations of the reflected Gaussian limit process: Weak convergence focuses on expectations, whereas large deviations characterize logarithmic tail asymptotics. In this work we partially fill this gap by providing a framework for large deviation asymptotics of multiclass feedforward queueing networks with first come first serve (fifo) and priority service disciplines at the nodes in critical loading. *Feedforward* means that the nodes can be numbered in such a way that the sequences of stations which are visited by customers have strictly increasing numbers. We consider deterministic class dependent routing of customers between nodes.

In the setting of single class queueing networks such a framework has been developed by Puhalskii [32]. There, a convergence toward critical loading is combined with a moderate deviation principle for the centered input processes. In this work a *moderate deviation principle* is regarded as a large deviation principle with a scaling close to a functional central limit theorem such that the rate function for the tail probabilities of a Gaussian process appears [6, 26, 35]. In this way Puhalskii obtains the characterization of logarithmic tail asymptotics for the appropriately scaled behavior of single class queueing networks when the queues approach a critical load. The rate function in the examples of Puhalskii [32] also appears in large deviation principles for reflected Brownian motion. Puhalskii's approach and its carryover to multiclass feedforward networks, presented here, is able to take advantage of both the simplified heavy traffic model captured through a Skorokhod map and the explicit rate function of a Gaussian process.

Puhalskii's work addresses random routing and feedback between the queues which is not permitted in our work where the feedforward structure allows for an iterative propagation of certain properties through the queues of the network. In our situation the difficulties come rather from the multiclass fifo service discipline. We note that in technical applications (e.g., the analysis of broadband communication networks) one is frequently faced with multiclass queueing networks where (at least in large parts of the network) feedback is not relevant, but the mixture of fifo and priority service disciplines.

The starting point for our results is a large deviation principle for the sequence of centered input processes which characterize the exogenous arrival and service times at the nodes. We then make use of the model for the behavior of multiclass feedforward queueing networks developed in [21]. There, it has been shown that close to a homogeneous critical load situation, this behavior can be approximated through a multi-dimensional Skorokhod map. We use this *pathwise heavy traffic convergence* to deduce sample path large deviation principles for the sequences of idle-time, workload, queue length, departure and sojourn time processes of the queueing network induced by the sequence of input processes. Besides the mentioned pathwise heavy traffic convergence, the proof of this result relies on a version of the contraction principle provided in [31].

In order to obtain these results, the mean arrival and service rates must be chosen in such a way that a critical load is approached at a certain speed at every node. In addition, our heavy traffic condition requires that the mean rates in the sequence of input processes grow to infinity. A moderate deviation principle can provide the basis for satisfying this heavy traffic condition, since it scales time faster than space. When our results are applied to input processes satisfying a moderate deviation principle, the rate function of a multi-dimensional reflected Gaussian process appears. This establishes a rigorous link between large deviations of the behavior of heavily loaded queueing networks and those of reflected Gaussian processes.

In this work we develop in detail this kind of relation for renewal arrival and service time processes. Our example partially parallels one of Puhalskii [32] and yields the rate function of a reflected Brownian motion which similarly appears in the weak heavy traffic approximation of Peterson [29]. In a companion work [22] we apply the results of this work to a two queue network with long-range dependent traffic. There, large and moderate deviation principles and explicit logarithmic tail asymptotics for the steady-state queue length at the second queue are obtained.

Our model captures the doubly infinite time interval. This permits us to treat directly steady-state situations. However, this feature requires the use of a topology which is slightly stronger than the usually considered topology of uniform convergence on compacts. In particular, the large deviation principle for the centered input processes must be provided in this topology before the main results of this work can be applied. By restricting the model to the positive time interval and by having initially empty queues, such a strengthening can be avoided. In addition, the considered model allows for mixtures of discrete customers and liquid inventory.

Our construction of the network behavior on the doubly infinite time interval works without a definition of the state of the network, which is not obvious for multiclass fifo queues containing discrete customers and fluid. Hence, we cannot deal with networks starting from a fixed or random nonzero state at time 0. By modifying the input processes on the negative time interval, one can force the network into certain nonzero initial states, but a detailed treatment of nonzero initial states goes beyond the scope of this work. We therefore refrain from claiming uniqueness of the stationary network behavior and refer to [5] for results in this direction.

An overview of this work is as follows: We recall necessary basic definitions and facts in Section 2. We state the main result for a single multiclass queue in Section 3. This result is extended to feedforward queueing networks in Section 4. We develop an application to renewal arrival and service time processes in Section 5. Appendix A establishes moderate deviation principles for renewal processes in the topology considered here. Appendix B presents conditions for stationarity of network behavior.

**2. Preliminaries.** A $[0,\infty]$-valued lower semicontinuous function on a topological Hausdorff space is called *rate function*. A rate function is *good* if it has compact level sets. A sequence $(X_k)_{k\in\mathbb{N}}$ of random elements with values in a measurable space $(E,\mathcal{E})$ satisfies a *large deviation principle* with normalizing sequence $(b_k)_{k\in\mathbb{N}}$ and rate function $I$ in the Hausdorff topology $\mathcal{T}$ on $E$ if $\lim_{k\to\infty} b_k = \infty$ and for every measurable set $A \in \mathcal{E}$,

$$\limsup_{k\to\infty} \frac{1}{b_k} \log P(X_k \in A) \leq -\inf_{x \in A^c} I(x)$$

and

$$\liminf_{k\to\infty} \frac{1}{b_k} \log P(X_k \in A) \geq -\inf_{x \in A^o} I(x),$$

where $A^c$ (resp. $A^o$) is the closure (resp. interior) of $A$ in the topology $\mathcal{T}$, and $P$ is the underlying probability measure. For background on large deviations theory, we refer to [7]. See also [12] for background on large deviations of queueing systems. In this work the random elements $X_k$ are continuous time processes and the space $E$ contains their sample paths which are introduced next.

We let $\mathcal{D}$ be the set of *functions* $\mathbf{d} : \mathbb{R} \to \mathbb{R}$ which are right continuous, possess left-hand limits and have limits

$$\underline{\mathbf{d}} := \lim_{t\to-\infty} \mathbf{d}(t)/t$$

and

$$\overline{\mathbf{d}} := \lim_{t\to\infty} \mathbf{d}(t)/t$$

in $\mathbb{R}$. On product function spaces these limits are taken componentwise. With $\mathcal{I}$ (resp. $\mathcal{C}$), we denote the subset of $\mathcal{D}$ consisting of nondecreasing (resp. continuous) functions. We provide these (and further) function spaces with the $\sigma$-algebras generated by the family of one-dimensional projections, and topology induced by the norm $\|\cdot\|$ given by

(1) $$\|\mathbf{d}\| := \sup_{t\in\mathbb{R}} \frac{|\mathbf{d}(t)|}{1+|t|}.$$

With the usual pointwise addition and scalar multiplication, the spaces $\mathcal{D}$ and $\mathcal{C}$ become topological vector spaces. Product spaces are equipped with corresponding product $\sigma$-algebras and topologies. We shall mostly use bold face symbols for elements of, or mappings into, function spaces.

We let $\mathcal{D}_0$ be the subset of $\mathcal{D}$ containing all functions $\mathbf{d} \in \mathcal{D}$ satisfying $\underline{\mathbf{d}} = \overline{\mathbf{d}} = 0$. Clearly, if a function $\mathbf{d} \in \mathcal{D}$ satisfies $\underline{\mathbf{d}} = \overline{\mathbf{d}} =: \delta$, then the *centered* function $\mathbf{d} - \delta \mathbf{id}$ is an element of $\mathcal{D}_0$, where $\mathbf{id} : \mathbb{R} \to \mathbb{R}, t \mapsto t$ denotes the identity map of $\mathbb{R}$.

We recall that the composition $\circ : (\mathcal{C} \times \mathcal{D}) \cup (\mathcal{D} \times \mathcal{I}) \to \mathcal{D}, (\mathbf{c}, \mathbf{d}) \mapsto \mathbf{c} \circ \mathbf{d}$ is continuous at elements in $\mathcal{C} \times \mathcal{D}$ and satisfies $\underline{\mathbf{c} \circ \mathbf{d}} = \underline{\mathbf{c}}\,\underline{\mathbf{d}}$ and $\overline{\mathbf{c} \circ \mathbf{d}} = \overline{\mathbf{c}}\,\overline{\mathbf{d}}$ (see Lemma B.1 in [19]).

On the set $\mathcal{D}_{\sup} := \{\mathbf{d} \in \mathcal{D} : \underline{\mathbf{d}} > 0\}$ we consider the "running supremum" $\sup : \mathcal{D}_{\sup} \to \mathcal{D}_{\sup} \cap \mathcal{I}$ defined for $\mathbf{d} \in \mathcal{D}_{\sup}$ and $t \in \mathbb{R}$ by

$$(\sup \mathbf{d})(t) := \sup_{\tau \in ]-\infty, t]} \mathbf{d}(\tau).$$

This mapping is continuous (see [11]), and satisfies $\underline{\sup \mathbf{d}} = \underline{\mathbf{d}}$ and $\overline{\sup \mathbf{d}} = \max\{0, \overline{\mathbf{d}}\}$. (For a finite subset $\mathcal{K}$ of $\mathbb{R}$, we let $\max \mathcal{K}$ denote the maximum of its elements.)

On the set $\mathcal{I}_{\mathrm{inv}} := \{\mathbf{c} \in \mathcal{I} \cap \mathcal{C} : \underline{\mathbf{c}} > 0, \overline{\mathbf{c}} > 0\}$ we define the "right continuous inverse" $\cdot^{-1} : \mathcal{I}_{\mathrm{inv}} \to \mathcal{I}$ for $\mathbf{c} \in \mathcal{I}_{\mathrm{inv}}$ and $t \in \mathbb{R}$ by

$$(\mathbf{c}^{-1})(t) := \sup_{\tau \in \mathbb{R}, \mathbf{c}(\tau) \leq t} \tau.$$

Clearly, $\mathbf{c} \in \mathcal{I}_{\mathrm{inv}}$ implies $\mathbf{c} \circ \mathbf{c}^{-1} = \mathbf{id} \leq \mathbf{c}^{-1} \circ \mathbf{c}$, $\underline{\mathbf{c}^{-1}} = 1/\underline{\mathbf{c}}$ and $\overline{\mathbf{c}^{-1}} = 1/\overline{\mathbf{c}}$. Here inequalities with vectors or functions are satisfied if they hold for every component or argument.

For sets $\mathcal{S}$ and $\mathcal{K}$, we let $\mathcal{S}^{\mathcal{K}}$ be a set of (column) vectors with indices in $\mathcal{K}$ and components in $\mathcal{S}$. In particular, $\mathcal{S}^k$ can be identified with $\mathcal{S}^{\{1,\ldots,k\}}$ for $k \in \mathbb{N}$. For $v \in \mathcal{S}^{\mathcal{K}}$ and $k \in \mathcal{K}$, we let $v_k \in \mathcal{S}$ be the component of $v$ with index $k \in \mathcal{K}$.

**3. Large deviations for a single multiclass queue.** In this section we model the behavior of a queueing node with infinite buffer capacity and multiple customer classes. We state a large deviation principle for the behavior of the queueing processes when the critical load 1 is approached, provided the input processes, which characterize arrival and service times, satisfy an appropriate sample path large deviation principle.

We let $\mathcal{M}$ be a finite set of customer classes. Two disjoint subsets $\mathcal{H}$ and $\mathcal{L}$ of the set $\mathcal{M}$ specify the high- and low-priority customer classes which are to be served by the node. Customers of a class in $\mathcal{M} \setminus (\mathcal{H} \cup \mathcal{L})$ do not have to visit the considered node and leave immediately. Such classes are allowed in order to facilitate the extension to feedforward networks in the following Section 4. Customers of a class in $\mathcal{H}$ (resp. $\mathcal{L}$) receive high-priority (resp. low-priority) service. Whenever a high-priority customer is present at the node, the service of low-priority customers is interrupted, and the customers with high-priority are served using a fifo service discipline, before service of low-priority customers is continued. Because of this *preemptive-resume* priority service discipline, the high-priority customers experience a service by the queueing node as if there were no low-priority customers. Customers with a class in $\mathcal{L}$ share the remaining service capacity of the queueing node using a fifo service discipline. We assume that $\mathcal{L}$ is nonempty. This still permits modeling a pure fifo service discipline (without priorities) at the node by letting $\mathcal{H}$ be empty.

The external arrivals at the node are specified through elements $\mathbf{a}$ of $\mathcal{I}^{\mathcal{M}}$. The value of the $j$th component of such a nondecreasing arrival vector function $\mathbf{a}$ at time $t$ is interpreted as the cumulative number of class $j$ arrivals up to time $t$. This means that $\mathbf{a}_j(t) - \mathbf{a}_j(s) \geq 0$ counts the number of arrivals of class $j$ customers during the time interval $]s, t]$ for $s < t \in \mathbb{R}$ and $j \in \mathcal{M}$. Here and in the following *cumulative* number of customers can be fractional (in order to permit scaling and to describe movements of fluid) and negative (in order to extend the cumulative customer numbers to the infinite negative past). Since the components of $\mathbf{a}$ are nondecreasing, there are no "negative" customers in the system. A jump in an arrival function $\mathbf{a}_j$ can be interpreted as the arrival of a discrete customer with a size corresponding to the height of the jump, whereas a continuous increase corresponds to a movement of fluid. The vectors $\underline{\mathbf{a}}, \overline{\mathbf{a}} \in \mathbb{R}_+^{\mathcal{M}}$ represent the asymptotical arrival rates on the negative and positive time intervals, respectively.

The service time requirements of the customers needing service at the node are determined through a vector function $\mathbf{s} \in \mathcal{I}_{\text{inv}}^{\mathcal{H} \cup \mathcal{L}}$. For a customer class $j \in \mathcal{H} \cup \mathcal{L}$, the value $\mathbf{s}_j(b)$ specifies the cumulative number of served customers of class $j$ when the node dedicates the total cumulative busy time $b \in \mathbb{R}$ uniquely to this class of customers. In particular, the function $\mathbf{s}_j^{-1}$ (resp. $\mathbf{s}_j^{-1} \circ \mathbf{a}_j$) describes the cumulative busy time required to serve the class $j$ customers as a function of the cumulative number of arrived class $j$ customers (resp. cumulative elapsed time). The continuity of $\mathbf{s}_j$ implies that the service time for a strictly positive number of class $j$ customers is strictly positive. Jumps in the function $\mathbf{s}_j^{-1}$ correspond to a service interruption for class $j$ customers. The vectors $\underline{\mathbf{s}}, \overline{\mathbf{s}} \in \mathbb{R}_+^{\mathcal{H} \cup \mathcal{L}}$ represent the asymptotical potential service rates on the negative and positive time intervals.

In the following queueing model the vector function $\mathbf{s}$ only appears in terms of the form $\mathbf{s}_j^{-1} \circ \mathbf{a}_j$ for $j \in \mathcal{M}$. Hence, only the values of $\mathbf{s}_j^{-1}$ on the image $\{\mathbf{a}_j(t) : t \in \mathbb{R}\}$ of $\mathbf{a}_j$ are relevant for the behavior of the system. When the function $\mathbf{a}_j$ for $j \in \mathcal{H} \cup \mathcal{L}$ is piecewise constant, the sequence of jump sizes of the function $\mathbf{s}_j^{-1} \circ \mathbf{a}_j$ can be interpreted as the service times required by the discrete customers modeled through $\mathbf{a}_j$. The behavior of the strictly increasing function $\mathbf{s}_j^{-1}$ between these jumps is irrelevant as long as the image of $\mathbf{a}_j$ does not change. In Section 5 we will use these facts to model the vector functions $\mathbf{a}$ and $\mathbf{s}$ from sequences of inter-arrival and service times.

Pairs $(\mathbf{a},\mathbf{s})\in\mathcal{I}^{\mathcal{M}}\times\mathcal{I}_{\mathrm{inv}}^{\mathcal{H}\cup\mathcal{L}}$ satisfying the *stability condition*

$$\sum_{j\in\mathcal{H}\cup\mathcal{L}}\frac{\underline{\mathbf{a}}_j}{\underline{\mathbf{s}}_j}<1 \tag{2}$$

are called *regular node primitives*. The stability condition prevents overload during the infinite past, but does not forbid overload during finite time intervals or the infinite future. Given a regular node primitives pair $(\mathbf{a},\mathbf{s})$, the cumulative idle time development of the node after service of only high-priority customers is defined as

$$\mathbf{U}^{\mathcal{H}}(\mathbf{a},\mathbf{s}):=\sup\left(\mathrm{id}-\sum_{j\in\mathcal{H}}\mathbf{s}_j^{-1}\circ\mathbf{a}_j\right)\in\mathcal{I}. \tag{3}$$

The (backlogged) workload behavior of high-priority customers (i.e., the time needed to complete service of all high-priority customers which are currently at the node) can be expressed as

$$\mathbf{V}^{\mathcal{H}}(\mathbf{a},\mathbf{s}):=\sum_{j\in\mathcal{H}}\mathbf{s}_j^{-1}\circ\mathbf{a}_j-\mathrm{id}+\mathbf{U}^{\mathcal{H}}(\mathbf{a},\mathbf{s})\in\mathcal{D}. \tag{4}$$

For $\mathcal{H}=\varnothing$ (pure fifo service discipline), we get $\mathbf{U}^{\mathcal{H}}(\mathbf{a},\mathbf{s})=\mathrm{id}$ and $\mathbf{V}^{\mathcal{H}}(\mathbf{a},\mathbf{s})=0$.

The remaining unused potential busy time of the node is dedicated to serving low-priority customers. Therefore, the cumulative idle time development of the node is defined by

$$\mathbf{Y}^{\mathcal{H},\mathcal{L}}(\mathbf{a},\mathbf{s}):=\sup\left(\mathbf{U}^{\mathcal{H}}(\mathbf{a},\mathbf{s})-\sum_{j\in\mathcal{L}}\mathbf{s}_j^{-1}\circ\mathbf{a}_j\right)\in\mathcal{I}. \tag{5}$$

The (backlogged) workload behavior of low-priority customers is given by

$$\mathbf{W}^{\mathcal{H},\mathcal{L}}(\mathbf{a},\mathbf{s}):=\sum_{j\in\mathcal{L}}\mathbf{s}_j^{-1}\circ\mathbf{a}_j-\mathbf{U}^{\mathcal{H}}(\mathbf{a},\mathbf{s})+\mathbf{Y}^{\mathcal{H},\mathcal{L}}(\mathbf{a},\mathbf{s})\in\mathcal{D}. \tag{6}$$

Lemma 2 of [21] shows that these definitions imply

$$\mathbf{V}^{\mathcal{H}}(\mathbf{a},\mathbf{s})\geq 0,$$
$$\underline{\mathbf{V}^{\mathcal{H}}(\mathbf{a},\mathbf{s})}=0,$$
$$\int_{\mathbb{R}}\mathbf{V}^{\mathcal{H}}(\mathbf{a},\mathbf{s})(t)\,d\mathbf{U}^{\mathcal{H}}(\mathbf{a},\mathbf{s})(t)=0,$$
$$\mathbf{W}^{\mathcal{H}}(\mathbf{a},\mathbf{s})\geq 0,$$
$$\underline{\mathbf{W}^{\mathcal{H}}(\mathbf{a},\mathbf{s})}=0,$$
$$\int_{\mathbb{R}}\mathbf{W}^{\mathcal{H},\mathcal{L}}(\mathbf{a},\mathbf{s})(t)\,d\mathbf{Y}^{\mathcal{H},\mathcal{L}}(\mathbf{a},\mathbf{s})(t)=0,$$
$$\int_{\mathbb{R}}\mathbf{V}^{\mathcal{H}}(\mathbf{a},\mathbf{s})(t)\,d\mathbf{Y}^{\mathcal{H},\mathcal{L}}(\mathbf{a},\mathbf{s})(t)=0.$$

Here $d\mathbf{u}$ denotes the $\sigma$-finite measure induced by a nondecreasing path $\mathbf{u} \in \mathcal{I}$ on $\mathbb{R}$. Conversely, this lemma also states that these properties in conjunction with equations (2), (4) and (6) imply equations (3) and (5). Hence, the definitions in equations (3) to (6) are justified by the nonidling service, the preemptive-resume priority structure of the service discipline, the wish for stability on the negative time interval, and the interpretation of $\mathbf{s}_j^{-1} \circ \mathbf{a}_j$ as the cumulative busy time required to serve all arrivals of class $j$ customers as a function of time.

Lemma 2 of [21] also establishes the hidden identities

$$\mathbf{Y}^{\mathcal{H},\mathcal{L}}(\mathbf{a},\mathbf{s}) = \sup\left(\mathbf{id} - \sum_{j \in \mathcal{H} \cup \mathcal{L}} \mathbf{s}_j^{-1} \circ \mathbf{a}_j\right),$$

$$\mathbf{W}^{\mathcal{H},\mathcal{L}}(\mathbf{a},\mathbf{s}) + \mathbf{V}^{\mathcal{H}}(\mathbf{a},\mathbf{s}) = \sum_{j \in \mathcal{H} \cup \mathcal{L}} \mathbf{s}_j^{-1} \circ \mathbf{a}_j - \mathbf{id} + \mathbf{Y}^{\mathcal{H},\mathcal{L}}(\mathbf{a},\mathbf{s}).$$

Next, we define the cumulative departures $\mathbf{D}^{\mathcal{H},\mathcal{L}}(\mathbf{a},\mathbf{s}) \in \mathcal{I}^{\mathcal{M}}$ at the node. The $j$th component in this vector of functions describes the cumulative number of customers of class $j \in \mathcal{M}$ which depart from the node as a function of time. We define it for $t \in \mathbb{R}$ by

$$\mathbf{D}_j^{\mathcal{H},\mathcal{L}}(\mathbf{a},\mathbf{s})(t)$$

$$:= \begin{cases} \sup_{\tau \leq t, \sum_{\ell \in \mathcal{H}} \mathbf{s}_\ell^{-1} \circ \mathbf{a}_\ell(\tau) \leq \sum_{\ell \in \mathcal{H}} \mathbf{s}_\ell^{-1} \circ \mathbf{a}_\ell(t) - \mathbf{V}^{\mathcal{H}}(\mathbf{a},\mathbf{s})(t)} \mathbf{a}_j(\tau), & \text{if } j \in \mathcal{H}, \\ \sup_{\tau \leq t, \sum_{\ell \in \mathcal{L}} \mathbf{s}_\ell^{-1} \circ \mathbf{a}_\ell(\tau) \leq \sum_{\ell \in \mathcal{L}} \mathbf{s}_\ell^{-1} \circ \mathbf{a}_\ell(t) - \mathbf{W}^{\mathcal{H},\mathcal{L}}(\mathbf{a},\mathbf{s})(t)} \mathbf{a}_j(\tau), & \text{if } j \in \mathcal{L}, \\ \mathbf{a}_j(t), & \text{otherwise.} \end{cases}$$

The conditions below the supremum in this definition assures that the workload incorporated in high- or low-priority customers at the node reaches at least its current workload. Due to partially served customers, it can exceed the current workload. Customers of each priority group leave in the order of their arrivals, and all customers of one priority group arriving at the same time instant leave together. Conversely, a larger or smaller number of departures would not be consistent with this order of the departures within the priority groups, the requirement that served customers must leave, or the definition of the workload behavior. In this way $\mathbf{D}^{\mathcal{H},\mathcal{L}}$ models fifo service discipline within the class of high- or low-priority customers. Customers of classes which need not visit the node leave immediately.

The behavior of queue lengths can now be defined as the difference between cumulative arrivals and departures

$$\mathbf{Q}^{\mathcal{H},\mathcal{L}}(\mathbf{a},\mathbf{s}) := \mathbf{a} - \mathbf{D}^{\mathcal{H},\mathcal{L}}(\mathbf{a},\mathbf{s}) \in \mathcal{D}^{\mathcal{M}}.$$

Finally, the time at which an *observer* (e.g., an additional customer which does not affect the behavior of the system) of class $j \in \mathcal{M}$ entering the node at time $t \in \mathbb{R}$ leaves the node is

$$\mathbf{Z}_j^{\mathcal{H},\mathcal{L}}(\mathbf{a},\mathbf{s})(t) := \begin{cases} \inf_{\tau \geq t,\ \mathbf{D}_\ell^{\mathcal{H},\mathcal{L}}(\tau) \geq \mathbf{a}_\ell(t) \text{ for every } \ell \in \mathcal{H}} \tau, & \text{if } j \in \mathcal{H}, \\ \inf_{\tau \geq t,\ \mathbf{D}_\ell^{\mathcal{H},\mathcal{L}}(\tau) \geq \mathbf{a}_\ell(t) \text{ for every } \ell \in \mathcal{L}} \tau, & \text{if } j \in \mathcal{L}, \\ t, & \text{otherwise.} \end{cases}$$

Hence, $\mathbf{Z}_j^{\mathcal{H},\mathcal{L}}(\mathbf{a},\mathbf{s}) - \mathbf{id} \in \mathcal{D}^\mathcal{M}$ represents the sojourn time of a class $j$ observer at the node as a function of its arrival time.

The model could be extended to incorporate more than two priorities; compare [4]. Since high-priority customers vanish in the heavy traffic limit (see Remark 4 after Theorem 1 below), we refrain from such an extension.

Similar pathwise models are frequently used to analyze the behavior of queueing systems [4, 13]. The one presented here is special in that it unifies discrete customer and fluid models. This is particularly useful in our context since sample path large deviation scalings often lead to a focus on continuous functions (i.e., the rate function is infinite for noncontinuous arguments) corresponding to movements of fluid through the system, even when the cumulative arrival (and, hence, also the departure) processes model discrete customers (i.e., they are concentrated on piecewise constant functions).

Our model includes the infinite negative time interval which is particularly advantageous in analyzing stationary behavior [11]. We refer to Lemma B.2 in Appendix B for independence and stationary increments conditions which are sufficient for a stationary behavior of the queue.

On the other hand, by restricting the elements of all considered function spaces to the time interval $\mathbb{R}_+$, and, furthermore, to those starting in 0, the model captures the behavior of the multiclass node when it starts empty at time 0. Clearly, in this modification, the supremum in the definitions of the mappings $\|\cdot\|$, $\mathbf{sup}$, and $\cdot^{-1}$ has to be restricted to nonnegative arguments, and condition (2) becomes meaningless.

The full model with time interval $\mathbb{R}$ and its restriction starting empty at time zero are intimately related: If $(\mathbf{a},\mathbf{s})$ is a regular node primitives pair and $t \in \mathbb{R}$ a time where the associated queue is empty [i.e., $\mathbf{W}^{\mathcal{H},\mathcal{L}}(\mathbf{a},\mathbf{s})(t) = \mathbf{V}^\mathcal{H}(\mathbf{a},\mathbf{s})(t) = 0$], then, for every $t' \geq t$, the supremum in the definitions of $\mathbf{U}^\mathcal{H}(\mathbf{a},\mathbf{s})(t')$, $\mathbf{Y}^{\mathcal{H},\mathcal{L}}(\mathbf{a},\mathbf{s})(t')$ and $\mathbf{D}^{\mathcal{H},\mathcal{L}}(\mathbf{a},\mathbf{s})(t')$ can be restricted to the finite interval $[t, t']$. This observation implies that if $(\hat{\mathbf{a}}, \hat{\mathbf{s}})$ is a random node primitives pair satisfying $(\hat{\mathbf{a}}(0), \hat{\mathbf{s}}(0)) = 0$ and the assumptions of Lemma B.2 (implying the stationarity of the associated behavior of the queue), workload and queue length distributions in the model restricted to the positive time interval, converge to the corresponding stationary distributions as time increases to infinity. This property carries over to the network case of Section 4.

Next, we specify a time-homogeneous critical load situation: For the remainder of this section, we fix vectors $\sigma \in \mathbb{R}_+^{\mathcal{H} \cup \mathcal{L}}$ and $\alpha \in \mathbb{R}_+^{\mathcal{M}}$ satisfying

$$\sigma > 0,$$

(7)
$$\rho^{\mathcal{L}} := \sum_{j \in \mathcal{L}} \frac{\alpha_j}{\sigma_j} > 0$$

and

(8)
$$\sum_{j \in \mathcal{H} \cup \mathcal{L}} \frac{\alpha_j}{\sigma_j} = 1.$$

If we interpret $\alpha$ as vector of arrival rates and $\sigma$ as vector of service rates, the last condition means that the node is critically loaded. Condition (7) is equivalent to the existence of a $j \in \mathcal{L}$ with $\alpha_j > 0$.

In order to formulate the main theorem of this section, we define for functions $\mathbf{a} \in \mathcal{D}^{\mathcal{M}}$ and $\mathbf{s} \in \mathcal{D}^{\mathcal{H} \cup \mathcal{L}}$ satisfying

(9)
$$\sum_{j \in \mathcal{H} \cup \mathcal{L}} \left( \frac{\mathbf{s}_j}{\sigma_j} \frac{\alpha_j}{\sigma_j} - \frac{\mathbf{a}_j}{\sigma_j} \right) > 0$$

the functions

$$\tilde{\mathbf{U}}^{\mathcal{H}}(\mathbf{a}, \mathbf{s}) := \sum_{j \in \mathcal{H}} \left( \frac{\mathbf{s}_j}{\sigma_j} \circ \left( \frac{\alpha_j}{\sigma_j} \mathbf{id} \right) - \frac{\mathbf{a}_j}{\sigma_j} \right) \in \mathcal{D},$$

$$\tilde{\mathbf{Y}}^{\mathcal{H},\mathcal{L}}(\mathbf{a}, \mathbf{s}) := \sup \left( \sum_{j \in \mathcal{H} \cup \mathcal{L}} \left( \frac{\mathbf{s}_j}{\sigma_j} \circ \left( \frac{\alpha_j}{\sigma_j} \mathbf{id} \right) - \frac{\mathbf{a}_j}{\sigma_j} \right) \right) \in \mathcal{I}$$

and

$$\tilde{\mathbf{W}}^{\mathcal{H},\mathcal{L}}(\mathbf{a}, \mathbf{s}) := \sum_{j \in \mathcal{H} \cup \mathcal{L}} \left( \frac{\mathbf{a}_j}{\sigma_j} - \frac{\mathbf{s}_j}{\sigma_j} \circ \left( \frac{\alpha_j}{\sigma_j} \mathbf{id} \right) \right) + \tilde{\mathbf{Y}}^{\mathcal{H},\mathcal{L}}(\mathbf{a}, \mathbf{s}) \in \mathcal{D}.$$

Furthermore, we let $\alpha^{\mathcal{L}} \in \mathbb{R}_+^{\mathcal{M}}$ be the vector with components

$$\alpha_j^{\mathcal{L}} := \begin{cases} \dfrac{\alpha_j}{\rho^{\mathcal{L}}}, & \text{if } j \in \mathcal{L}, \\ 0, & \text{if } j \in \mathcal{M} \setminus \mathcal{L}, \end{cases}$$

and $e^{\mathcal{L}}$ the vector in $\mathbb{R}^{\mathcal{M}}$ with $e_j^{\mathcal{L}} = 1$, if $j \in \mathcal{L}$, and $e_j^{\mathcal{L}} = 0$, otherwise.

THEOREM 1. *We let $(\hat{\mathbf{a}}_k)_{k \in \mathbb{N}}$ be a sequence of processes on $\mathcal{I}^{\mathcal{M}}$ and $(\hat{\mathbf{s}}_k)_{k \in \mathbb{N}}$ a sequence of processes on $\mathcal{I}_{\text{inv}}^{\mathcal{H} \cup \mathcal{L}}$. We assume that, for every $k \in \mathbb{N}$, there are vectors $\hat{\alpha}_k \in \mathbb{R}_+^{\mathcal{M}}$ and $\hat{\sigma}_k \in \mathbb{R}_+^{\mathcal{H} \cup \mathcal{L}}$ such that, with probability one,*

$$\hat{\alpha}_k = \underline{\hat{\mathbf{a}}}_k = \overline{\hat{\mathbf{a}}}_k,$$

$$\hat{\sigma}_k = \underline{\hat{\mathbf{s}}}_k = \overline{\hat{\mathbf{s}}}_k > 0,$$

(10)
$$\hat{\rho}_k := \sum_{j \in \mathcal{H} \cup \mathcal{L}} \frac{\hat{\alpha}_{k,j}}{\hat{\sigma}_{k,j}} < 1.$$

We assume that the sequence of centered regular node primitives processes

$$(\hat{\mathbf{a}}_k - \hat{\alpha}_k \mathbf{id}, \hat{\mathbf{s}}_k - \hat{\sigma}_k \mathbf{id})_{k \in \mathbb{N}}$$

satisfies a sample path large deviation principle with normalizing sequence $(b_k)_{k \in \mathbb{N}}$ and good rate function $I$ on the space $\mathcal{D}_0^{\mathcal{M}} \times \mathcal{D}_0^{\mathcal{H} \cup \mathcal{L}}$. We assume that $I$ takes the value $\infty$ whenever a component of its argument functions is not continuous.

If there is a sequence $(d_k)_{k \in \mathbb{N}}$ in $\mathbb{R}_+$ and vectors $\tilde{\alpha} \in \mathbb{R}^{\mathcal{M}}$, $\tilde{\sigma} \in \mathbb{R}^{\mathcal{H} \cup \mathcal{L}}$ with

(11) $$\infty = \lim_{k \to \infty} d_k,$$

(12) $$\tilde{\alpha} = \lim_{k \to \infty} (\hat{\alpha}_k - d_k \alpha t),$$

(13) $$\tilde{\sigma} = \lim_{k \to \infty} (\hat{\sigma}_k - d_k \sigma)$$

and

(14) $$\tilde{\rho} := \sum_{j \in \mathcal{H} \cup \mathcal{L}} \left( \frac{\tilde{\sigma}_j}{\sigma_j} \frac{\alpha_j}{\sigma_j} - \frac{\tilde{\alpha}_j}{\sigma_j} \right) > 0,$$

the following 8 statements are valid with normalizing sequence $(b_k)_{k \in \mathbb{N}}$:

1. The sequence $(d_k \mathbf{U}^{\mathcal{H}}(\hat{\mathbf{a}}_k, \hat{\mathbf{s}}_k) - d_k(1 - \hat{\rho}_k^{\mathcal{H}})\mathbf{id})_{k \in \mathbb{N}}$, where

$$\hat{\rho}_k^{\mathcal{H}} := \sum_{j \in \mathcal{H}} \frac{\hat{\alpha}_{k,j}}{\hat{\sigma}_{k,j}},$$

satisfies a large deviation principle on $\mathcal{D}_0$ with good rate function

$$\tilde{\mathbf{u}} \mapsto \inf_{\tilde{\mathbf{a}} \in \mathcal{D}_0^{\mathcal{M}}, \tilde{\mathbf{s}} \in \mathcal{D}_0^{\mathcal{H} \cup \mathcal{L}},\ \tilde{\mathbf{u}} = \tilde{\mathbf{U}}^{\mathcal{H}}(\tilde{\mathbf{a}} + \tilde{\alpha}\mathbf{id}, \tilde{\mathbf{s}} + \tilde{\sigma}\mathbf{id}) - \tilde{\rho}^{\mathcal{H}}\mathbf{id}} I(\tilde{\mathbf{a}}, \tilde{\mathbf{s}}),$$

where

$$\tilde{\rho}^{\mathcal{H}} := \sum_{j \in \mathcal{H}} \left( \frac{\tilde{\sigma}_j}{\sigma_j} \frac{\alpha_j}{\sigma_j} - \frac{\tilde{\alpha}_j}{\sigma_j} \right).$$

2. The sequence $(d_k \mathbf{V}^{\mathcal{H}}(\hat{\mathbf{a}}_k, \hat{\mathbf{s}}_k))_{k \in \mathbb{N}}$ satisfies a large deviation principle on $\mathcal{D}_0$ with good rate function

$$\tilde{\mathbf{v}} \mapsto \begin{cases} 0, & \text{if } \tilde{\mathbf{v}} = 0, \\ \infty, & \text{otherwise.} \end{cases}$$

3. The sequence $(d_k \mathbf{Y}^{\mathcal{H},\mathcal{L}}(\hat{\mathbf{a}}_k, \hat{\mathbf{s}}_k) - d_k(1 - \hat{\rho}_k)\mathbf{id})_{k \in \mathbb{N}}$ satisfies a large deviation principle on $\mathcal{D}_0$ with good rate function

$$\tilde{\mathbf{y}} \mapsto \inf_{\tilde{\mathbf{a}} \in \mathcal{D}_0^{\mathcal{M}}, \tilde{\mathbf{s}} \in \mathcal{D}_0^{\mathcal{H} \cup \mathcal{L}}, \tilde{\mathbf{y}} = \tilde{\mathbf{Y}}^{\mathcal{H},\mathcal{L}}(\tilde{\mathbf{a}} + \tilde{\alpha}\mathbf{id}, \tilde{\mathbf{s}} + \tilde{\sigma}\mathbf{id}) - \tilde{\rho}\mathbf{id}} I(\tilde{\mathbf{a}}, \tilde{\mathbf{s}}).$$

4. *The sequence $(d_k \mathbf{W}^{\mathcal{H},\mathcal{L}}(\hat{\mathbf{a}}_k, \hat{\mathbf{s}}_k))_{k\in\mathbb{N}}$ satisfies a large deviation principle on $\mathcal{D}_0$ with good rate function*

$$\tilde{\mathbf{w}} \mapsto \inf_{\tilde{\mathbf{a}}\in\mathcal{D}_0^{\mathcal{M}},\tilde{\mathbf{s}}\in\mathcal{D}_0^{\mathcal{H}\cup\mathcal{L}},\tilde{\mathbf{w}}=\tilde{\mathbf{W}}^{\mathcal{H},\mathcal{L}}(\tilde{\mathbf{a}}+\tilde{\alpha}\mathbf{id},\tilde{\mathbf{s}}+\tilde{\sigma}\mathbf{id})} I(\tilde{\mathbf{a}},\tilde{\mathbf{s}}).$$

5. *The sequence $(\mathbf{D}^{\mathcal{H},\mathcal{L}}(\hat{\mathbf{a}}_k, \hat{\mathbf{s}}_k) - \hat{\alpha}_k \mathbf{id})_{k\in\mathbb{N}}$ satisfies a large deviation principle on $\mathcal{D}_0^{\mathcal{M}}$ with good rate function*

$$\tilde{\mathbf{d}} \mapsto \inf_{\tilde{\mathbf{a}}\in\mathcal{D}_0^{\mathcal{M}},\tilde{\mathbf{s}}\in\mathcal{D}_0^{\mathcal{H}\cup\mathcal{L}},\tilde{\mathbf{d}}=\tilde{\mathbf{a}}-\alpha^{\mathcal{L}}\tilde{\mathbf{W}}^{\mathcal{H},\mathcal{L}}(\tilde{\mathbf{a}}+\tilde{\alpha}\mathbf{id},\tilde{\mathbf{s}}+\tilde{\sigma}\mathbf{id})} I(\tilde{\mathbf{a}},\tilde{\mathbf{s}}).$$

6. *The sequence $(\mathbf{Q}^{\mathcal{H},\mathcal{L}}(\hat{\mathbf{a}}_k, \hat{\mathbf{s}}_k))_{k\in\mathbb{N}}$ satisfies a large deviation principle on $\mathcal{D}_0^{\mathcal{M}}$ with good rate function*

$$\tilde{\mathbf{q}} \mapsto \inf_{\tilde{\mathbf{a}}\in\mathcal{D}_0^{\mathcal{M}},\tilde{\mathbf{s}}\in\mathcal{D}_0^{\mathcal{H}\cup\mathcal{L}},\tilde{\mathbf{q}}=\alpha^{\mathcal{L}}\tilde{\mathbf{W}}^{\mathcal{H},\mathcal{L}}(\tilde{\mathbf{a}}+\tilde{\alpha}\mathbf{id},\tilde{\mathbf{s}}+\tilde{\sigma}\mathbf{id})} I(\tilde{\mathbf{a}},\tilde{\mathbf{s}}).$$

7. *The sequence $(d_k \mathbf{Z}^{\mathcal{H},\mathcal{L}}(\hat{\mathbf{a}}_k, \hat{\mathbf{s}}_k) - d_k \mathbf{id})_{k\in\mathbb{N}}$ satisfies a large deviation principle on $\mathcal{D}_0^{\mathcal{M}}$ with good rate function*

$$\tilde{\mathbf{z}} \mapsto \inf_{\tilde{\mathbf{a}}\in\mathcal{D}_0^{\mathcal{M}},\tilde{\mathbf{s}}\in\mathcal{D}_0^{\mathcal{H}\cup\mathcal{L}},\tilde{\mathbf{z}}=e^{\mathcal{L}}\tilde{\mathbf{W}}^{\mathcal{H},\mathcal{L}}(\tilde{\mathbf{a}}+\tilde{\alpha}\mathbf{id},\tilde{\mathbf{s}}+\tilde{\sigma}\mathbf{id})} I(\tilde{\mathbf{a}},\tilde{\mathbf{s}}).$$

8. *The previous statements are also valid in combination: For instance, the sequence $(d_k \mathbf{W}^{\mathcal{H},\mathcal{L}}(\hat{\mathbf{a}}_k, \hat{\mathbf{s}}_k), \mathbf{Q}^{\mathcal{H},\mathcal{L}}(\hat{\mathbf{a}}_k, \hat{\mathbf{s}}_k))_{k\in\mathbb{N}}$ satisfies a large deviation principle on $\mathcal{D}_0 \times \mathcal{D}_0^{\mathcal{M}}$ with good rate function*

$$(\tilde{\mathbf{w}}, \tilde{\mathbf{q}}) \mapsto \inf_{\substack{\tilde{\mathbf{a}}\in\mathcal{D}_0^{\mathcal{M}},\tilde{\mathbf{s}}\in\mathcal{D}_0^{\mathcal{H}\cup\mathcal{L}},\tilde{\mathbf{w}}=\tilde{\mathbf{W}}^{\mathcal{H},\mathcal{L}}(\tilde{\mathbf{a}}+\tilde{\alpha}\mathbf{id},\tilde{\mathbf{s}}+\tilde{\sigma}\mathbf{id}),\\ \tilde{\mathbf{q}}=\alpha^{\mathcal{L}}\tilde{\mathbf{W}}^{\mathcal{H},\mathcal{L}}(\tilde{\mathbf{a}}+\tilde{\alpha}\mathbf{id},\tilde{\mathbf{s}}+\tilde{\sigma}\mathbf{id})}} I(\tilde{\mathbf{a}},\tilde{\mathbf{s}}).$$

*The statements of the theorem remain valid when the elements in all underlying function spaces are restricted to the time interval $\mathbb{R}_+$ and to those starting in 0, and the definitions of the model are modified accordingly. [Hence, conditions (2) and (9) become meaningless.] With this modification, the statements are also true when conditions (10) and (14) are dropped, and/or the topology is changed to the topology of uniform convergence on compacts.*

Before giving the proof we remark:

REMARK 1. Limits (11) to (13) imply

$$\lim_{k\to\infty} \frac{\hat{\alpha}_k}{d_k} = \alpha,$$

$$\lim_{k\to\infty} \frac{\hat{\sigma}_k}{d_k} = \sigma,$$

$$\lim_{k\to\infty} \hat{\rho}_k = 1$$

and

$$\lim_{k\to\infty} d_k(1-\hat{\rho}_k) = \lim_{k\to\infty} \sum_{j\in\mathcal{H}\cup\mathcal{L}} \frac{(\hat{\sigma}_{k,j} - d_k\sigma_j)\alpha_j - (\hat{\alpha}_{k,j} - d_k\alpha_j)\sigma_j}{\sigma_j\hat{\sigma}_{k,j}/d_k} = \tilde{\rho}.$$

Hence, the asymptotic load at the queue on large time intervals must approach the critical value 1 with a speed of order $1/d_k$ as $k \to \infty$. This prerequisite is the reason for speaking about large deviations "in critical loading" as in [32].

REMARK 2. The requirement (14) is not as restrictive as it might appear since conditions (10) to (13) imply $\tilde{\rho} \geq 0$.

REMARK 3. The mappings $\tilde{\mathbf{U}}^{\mathcal{H}}$, $\tilde{\mathbf{Y}}^{\mathcal{H},\mathcal{L}}$ and $\tilde{\mathbf{W}}^{\mathcal{H},\mathcal{L}}$ appearing in the expressions of the rate functions for the network behavior have a much simpler structure than the network mappings (particularly the mapping $\mathbf{D}^{\mathcal{H},\mathcal{L}}$ incorporating the queueing discipline). They consist mainly of linear combinations of the argument functions and the **sup**-mapping. This reflects the fact that the behavior of the queueing node simplifies close to a time-homogeneous critical load; compare [21].

REMARK 4. In statement 8 of the theorem the rate function for workload and queue length behavior is only finite when the $i$th component of the queue length behavior is $\alpha_i^{\mathcal{L}}$ times the workload behavior. In this case, high-priority customers do not experience a delay. An analogous *state space collapse* is well known from heavy traffic convergence in distribution [29].

PROOF. We first prove statement 4. Under the assumptions of the theorem, we define for every $k \in \mathbb{N}$ the mapping $\mathbf{F}_k \colon \mathcal{D}_0^{\mathcal{M}} \times \mathcal{D}_0^{\mathcal{H}\cup\mathcal{L}} \to \mathcal{D}_0$ by

$$\mathbf{F}_k(\tilde{\mathbf{a}}, \tilde{\mathbf{s}}) := d_k \mathbf{W}^{\mathcal{H},\mathcal{L}}\left(\frac{\tilde{\mathbf{a}} + \hat{\alpha}_k \mathbf{id}}{d_k}, \frac{\tilde{\mathbf{s}} + \hat{\sigma}_k \mathbf{id}}{d_k}\right),$$

if $\tilde{\mathbf{a}} \in \mathcal{D}_0^{\mathcal{M}}$ and $\tilde{\mathbf{s}} \in \mathcal{D}_0^{\mathcal{H}\cup\mathcal{L}}$ are such that $\tilde{\mathbf{a}} + \hat{\alpha}_k \mathbf{id} \in \mathcal{I}^{\mathcal{M}}$ and $\tilde{\mathbf{s}} + \hat{\sigma}_k \mathbf{id} \in \mathcal{I}_{\mathrm{inv}}^{\mathcal{H}\cup\mathcal{L}}$. Otherwise, we set

$$\mathbf{F}_k(\tilde{\mathbf{a}}, \tilde{\mathbf{s}}) := \tilde{\mathbf{W}}^{\mathcal{H},\mathcal{L}}(\tilde{\mathbf{a}} + \tilde{\alpha}\mathbf{id}, \tilde{\mathbf{s}} + \tilde{\sigma}\mathbf{id}).$$

Conditions (10) and (14) ensure that the mapping $\mathbf{F}_k$ is well defined.

We let $\tilde{\mathbf{a}} \in \mathcal{D}_0^{\mathcal{M}}$ and $\tilde{\mathbf{s}} \in \mathcal{D}_0^{\mathcal{H}\cup\mathcal{L}}$ with $I(\tilde{\mathbf{a}}, \tilde{\mathbf{s}}) < \infty$. Hence, $\tilde{\mathbf{a}} \in \mathcal{C}^{\mathcal{M}}$ and $\tilde{\mathbf{s}} \in \mathcal{C}^{\mathcal{H}\cup\mathcal{L}}$ because $I(\tilde{\mathbf{a}}, \tilde{\mathbf{s}})$ would be infinite if a component of $\tilde{\mathbf{a}}$ or $\tilde{\mathbf{s}}$ were not continuous. If $(\tilde{\mathbf{a}}_k)_{k\in\mathbb{N}}$ is a sequence of functions in $\mathcal{D}_0^{\mathcal{M}}$ converging to $\tilde{\mathbf{a}}$ and $(\tilde{\mathbf{s}}_k)_{k\in\mathbb{N}}$ is a sequence of functions in $\mathcal{D}_0^{\mathcal{H}\cup\mathcal{L}}$ converging to $\tilde{\mathbf{s}}$, then

$$(15) \qquad \lim_{k\to\infty} \mathbf{F}_k(\tilde{\mathbf{a}}_k, \tilde{\mathbf{s}}_k) = \tilde{\mathbf{W}}^{\mathcal{H},\mathcal{L}}(\tilde{\mathbf{a}} + \tilde{\alpha}\mathbf{id}, \tilde{\mathbf{s}} + \tilde{\sigma}\mathbf{id}).$$

This is true if $\tilde{\mathbf{a}}_k + \hat{\alpha}_k \mathbf{id} \in \mathcal{I}^{\mathcal{M}}$ and $\tilde{\mathbf{s}}_k + \hat{\sigma}_k \mathbf{id} \in \mathcal{I}_{\text{inv}}^{\mathcal{H} \cup \mathcal{L}}$ for every $k \in \mathbb{N}$ because of the following: In this case we set, for $k \in \mathbb{N}$,

$$\mathbf{a}_k := \frac{\tilde{\mathbf{a}}_k + \hat{\alpha}_k \mathbf{id}}{d_k}$$

and

$$\mathbf{s}_k := \frac{\tilde{\mathbf{s}}_k + \hat{\sigma}_k \mathbf{id}}{d_k},$$

and thus obtain a sequence of regular node primitives $(\mathbf{a}_k, \mathbf{s}_k)_{k \in \mathbb{N}}$ which satisfies

$$\lim_{k \to \infty} d_k(\mathbf{a}_k - \alpha \mathbf{id}) = \lim_{k \to \infty} (\tilde{\mathbf{a}}_k + (\hat{\alpha}_k - d_k \alpha) \mathbf{id}) = \tilde{\mathbf{a}} + \tilde{\alpha} \mathbf{id}$$

and

$$\lim_{k \to \infty} d_k(\mathbf{s}_k - \sigma \mathbf{id}) = \lim_{k \to \infty} (\tilde{\mathbf{s}}_k + (\hat{\alpha}_k - d_k \sigma) \mathbf{id}) = \tilde{\mathbf{s}} + \tilde{\sigma} \mathbf{id},$$

because of conditions (12) and (13). Thus, this sequence satisfies the prerequisites of Theorem 4 in [21] and we obtain

$$\lim_{k \to \infty} \mathbf{F}_k(\tilde{\mathbf{a}}_k, \tilde{\mathbf{s}}_k) = \lim_{k \to \infty} d_k \mathbf{W}^{\mathcal{H}, \mathcal{L}}(\mathbf{a}_k, \mathbf{s}_k) = \tilde{\mathbf{W}}^{\mathcal{H}, \mathcal{L}}(\tilde{\mathbf{a}} + \tilde{\alpha} \mathbf{id}, \tilde{\mathbf{s}} + \tilde{\sigma} \mathbf{id}).$$

If $\tilde{\mathbf{a}}_k + \hat{\alpha}_k \notin \mathcal{I}^{\mathcal{M}}$ or $\tilde{\mathbf{s}}_k + \hat{\sigma}_k \notin \mathcal{I}_{\text{inv}}^{\mathcal{H} \cup \mathcal{L}}$ for some $k \in \mathbb{N}$, we still obtain the convergence (15) since $\tilde{\mathbf{W}}^{\mathcal{H}, \mathcal{L}}$ is continuous.

We can therefore apply the version of the contraction principle of [31] (see Remarks 2 and 3 after Theorem 2.1 therein) in order to conclude that the sequence $\mathbf{F}_k(\hat{\mathbf{a}}_k - \hat{\alpha}_k \mathbf{id}, \hat{\mathbf{s}}_k - \hat{\sigma}_k \mathbf{id})$ satisfies a sample path large deviation principle with the normalizing sequence $(b_k)_{k \in \mathbb{N}}$ and the good rate function given in statement 4. Statement 4 now follows from the fact that, for every $k \in \mathbb{N}$, we have, with probability one,

$$\mathbf{F}_k(\hat{\mathbf{a}}_k - \hat{\alpha}_k \mathbf{id}, \hat{\mathbf{s}}_k - \hat{\sigma}_k \mathbf{id}) = d_k \mathbf{W}^{\mathcal{H}, \mathcal{L}}\left(\frac{\hat{\mathbf{a}}_k}{d_k}, \frac{\hat{\mathbf{s}}_k}{d_k}\right) = d_k \mathbf{W}^{\mathcal{H}, \mathcal{L}}(\hat{\mathbf{a}}_k, \hat{\mathbf{s}}_k).$$

In the last equation we used the property of the mapping $\mathbf{W}^{\mathcal{H}, \mathcal{L}}$ that for every regular node primitives tuple $(\mathbf{a}, \mathbf{s})$ and value $\xi > 0$, we have

$$\mathbf{W}^{\mathcal{H}, \mathcal{L}}(\xi \mathbf{a}, \xi \mathbf{s}) = \mathbf{W}^{\mathcal{H}, \mathcal{L}}(\mathbf{a}, \mathbf{s}).$$

The proofs of the other large deviation principles are analogous. However, in certain cases the previous invariance must be replaced by the positive homogeneity as exemplified by the equation

$$\mathbf{D}^{\mathcal{H}, \mathcal{L}}(\xi \mathbf{a}, \xi \mathbf{s}) = \xi \, \mathbf{D}^{\mathcal{H}, \mathcal{L}}(\mathbf{a}, \mathbf{s}).$$

This distinction is responsible for the fact that in statements 5, 6 and 8 the factor $d_k$ disappears.

When the model is restricted to the time interval $\mathbb{R}_+$ and to functions starting at 0, the proofs are completely analog and do not make use of conditions (10) and (14). Compare Remarks 5 and 6 after Theorem 1 in [21]. In this situation the behavior of the queue on the time interval $[0,t]$ for $t \geq 0$ depends only on the cumulative arrivals until time $t$ and the cumulative number of served customers up to the cumulative busy time $t$. This permits us to transfer the statements to the topology of uniform convergence on compacts and completes the proof. $\square$

**4. Large deviations for feedforward networks.** In this section we consider $n \in \mathbb{N}$ queues numbered from 1 to $n$ with fifo and priority service disciplines. Therefore, several symbols from the previous section receive an additional index $i$ specifying the queue with which they are associated. Each customer of the network belongs to a class of the finite set of customer classes $\mathcal{M}$. For every queue $i \in \{1, \ldots, n\}$, disjoint subsets $\mathcal{H}_i$ and $\mathcal{L}_i$ of $\mathcal{M}$ specify the classes of customers which must visit the $i$th queue. Before leaving the network, the customers of class $j \in \mathcal{M}$ have to visit all queues $i \leq n$ with $j \in \mathcal{H}_i \cup \mathcal{L}_i$ in increasing order. Because of this last requirement, our model only captures feedforward networks.

Customers of a class in $\mathcal{H}_i$ receive high-priority and customers of a class in $\mathcal{L}_i$ low-priority service at queue $i \leq n$. In particular, the priority group to which a customer class belongs can change from station to station. When high-priority customers are present at a queue, the service of low-priority customers is interrupted and then resumed when no high-priority customers require service. Among the customers of one priority group the service time is shared in fifo order. We assume that $\mathcal{L}_i$ is nonempty for every $i \leq n$. By having $\mathcal{H}_i = \varnothing$ for some or all nodes $i \leq n$, one can model stations or networks with a pure fifo service discipline.

As for the single node, external arrivals at the network are specified through elements $\mathbf{a}$ of $\mathcal{I}^{\mathcal{M}}$. The service time requirements of the customers at the queues they have to visit are determined through vectors of functions $\mathbf{s}_1 \in \mathcal{I}_{\text{inv}}^{\mathcal{H}_1 \cup \mathcal{L}_1}, \ldots, \mathbf{s}_n \in \mathcal{I}_{\text{inv}}^{\mathcal{H}_n \cup \mathcal{L}_n}$. For every queue $i \leq n$ and customer class $j \in \mathcal{H}_i \cup \mathcal{L}_i$, the value of the component $\mathbf{s}_{i,j}$ at $b$ specifies the cumulative number of class $j$ customers which is served by queue $i$ when it dedicates the cumulative busy time $b$ to this class of customers.

In order to be able to construct the behavior of the queueing network for given cumulative arrival function vector $\mathbf{a}$ and service time function vectors $\mathbf{s}_1, \ldots, \mathbf{s}_n$, we have to guarantee that the queues have not been overloaded during the infinite past. This amounts to requiring the stability condition

$$\sum_{j \in \mathcal{H}_i \cup \mathcal{L}_i} \frac{\underline{\mathbf{a}}_j}{\underline{\mathbf{s}}_{i,j}} < 1 \qquad \text{for } i = 1, \ldots, n. \tag{16}$$

An $(n+1)$-tuple of function vectors $(\mathbf{a}, \mathbf{s}_1, \ldots, \mathbf{s}_n)$ with $\mathbf{a} \in \mathcal{I}^{\mathcal{M}}$, $\mathbf{s}_1 \in \mathcal{I}_{\text{inv}}^{\mathcal{H}_1 \cup \mathcal{L}_1}$, $\ldots$, $\mathbf{s}_n \in \mathcal{I}_{\text{inv}}^{\mathcal{H}_n \cup \mathcal{L}_n}$ satisfying this condition is called a *regular network primitives tuple*. For such a regular network primitives tuple $(\mathbf{a}, \mathbf{s}_1, \ldots, \mathbf{s}_n)$, we define

$$\mathbf{D}_0(\mathbf{a}, \mathbf{s}_1, \ldots, \mathbf{s}_n) := \mathbf{a}$$

and, inductively, for $i = 1, \ldots, n$,

$$\mathbf{U}_i(\mathbf{a}, \mathbf{s}_1, \ldots, \mathbf{s}_n) := \mathbf{U}^{\mathcal{H}_i}(\mathbf{D}_{i-1}(\mathbf{a}, \mathbf{s}_1, \ldots, \mathbf{s}_n), \mathbf{s}_i),$$
$$\mathbf{V}_i(\mathbf{a}, \mathbf{s}_1, \ldots, \mathbf{s}_n) := \mathbf{V}^{\mathcal{H}_i}(\mathbf{D}_{i-1}(\mathbf{a}, \mathbf{s}_1, \ldots, \mathbf{s}_n), \mathbf{s}_i),$$
$$\mathbf{Y}_i(\mathbf{a}, \mathbf{s}_1, \ldots, \mathbf{s}_n) := \mathbf{Y}^{\mathcal{H}_i, \mathcal{L}_i}(\mathbf{D}_{i-1}(\mathbf{a}, \mathbf{s}_1, \ldots, \mathbf{s}_n), \mathbf{s}_i),$$
$$\mathbf{W}_i(\mathbf{a}, \mathbf{s}_1, \ldots, \mathbf{s}_n) := \mathbf{W}^{\mathcal{H}_i, \mathcal{L}_i}(\mathbf{D}_{i-1}(\mathbf{a}, \mathbf{s}_1, \ldots, \mathbf{s}_n), \mathbf{s}_i),$$
$$\mathbf{D}_i(\mathbf{a}, \mathbf{s}_1, \ldots, \mathbf{s}_n) := \mathbf{D}^{\mathcal{H}_i, \mathcal{L}_i}(\mathbf{D}_{i-1}(\mathbf{a}, \mathbf{s}_1, \ldots, \mathbf{s}_n), \mathbf{s}_i),$$
$$\mathbf{Q}_i(\mathbf{a}, \mathbf{s}_1, \ldots, \mathbf{s}_n) := \mathbf{Q}^{\mathcal{H}_i, \mathcal{L}_i}(\mathbf{D}_{i-1}(\mathbf{a}, \mathbf{s}_1, \ldots, \mathbf{s}_n), \mathbf{s}_i)$$

and

$$\mathbf{Z}_i(\mathbf{a}, \mathbf{s}_1, \ldots, \mathbf{s}_n) := \mathbf{Z}^{\mathcal{H}_i, \mathcal{L}_i}(\mathbf{D}_{i-1}(\mathbf{a}, \mathbf{s}_1, \ldots, \mathbf{s}_n), \mathbf{s}_i).$$

These mappings specify the behavior of the $i$th node in the network analogously as in the single queue case. Again the model contains stationary (see Appendix B) and initially empty queueing networks (by restricting the functions in the underlying function spaces to the time interval $\mathbb{R}_+$ and to those starting in 0, and modifying the definitions of the model accordingly) as special cases. In the latter case condition (16) becomes meaningless.

We let vectors $\alpha \in \mathbb{R}_+^{\mathcal{M}}$, $\sigma_1 \in \mathbb{R}_+^{\mathcal{H}_1 \cup \mathcal{L}_1}$, $\ldots$, $\sigma_n \in \mathbb{R}_+^{\mathcal{H}_n \cup \mathcal{L}_n}$ be given which satisfy, for every $i \leq n$,

$$\sigma_i > 0,$$

$$(17) \qquad \rho_i^{\mathcal{L}} := \sum_{j \in \mathcal{L}_i} \frac{\alpha_j}{\sigma_{i,j}} > 0$$

and

$$(18) \qquad \sum_{j \in \mathcal{H}_i \cup \mathcal{L}_i} \frac{\alpha_j}{\sigma_{i,j}} = 1.$$

When $\alpha_j$ is interpreted as the exogenous arrival rate of class $j \in \mathcal{H}_i \cup \mathcal{L}_i$ customers and $\sigma_{i,j}$ as the service rate of class $j$ customers at queue $i \leq n$, the last condition says that node $i$ is critically loaded.

We shall need the vector $\alpha_i^{\mathcal{L}} \in \mathbb{R}^{\mathcal{M}}$ for $i \leq n$, which is defined by

$$(19) \qquad \alpha_{i,j}^{\mathcal{L}} := \begin{cases} \dfrac{\alpha_j}{\rho_i^{\mathcal{L}}}, & \text{if } j \in \mathcal{L}_i, \\ 0, & \text{if } j \in \mathcal{M} \setminus \mathcal{L}_i. \end{cases}$$

Furthermore, we define the strictly lower triangular matrix $G \in \mathbb{R}^{n \times n}$ by

$$(20) \quad G_{i,h} := \begin{cases} \sum_{j \in \mathcal{H}_i \cup \mathcal{L}_i} \dfrac{\alpha^{\mathcal{L}}_{h,j}}{\sigma_{i,j}}, & \text{if } h < i, \\ 0, & \text{otherwise.} \end{cases}$$

The triangularity of the matrix $G$ is a consequence of the feedforward structure of the modeled network.

Close to the "critical" linear network primitives tuple $(\alpha\mathbf{id}, \sigma_1\mathbf{id}, \ldots, \sigma_n\mathbf{id})$, the behavior of the network simplifies [21]. In order to model the simplified behavior we introduce for $i = 1, \ldots, n$ the $\mathcal{D}$-valued mappings $\tilde{\mathbf{X}}_i$, $\tilde{\mathbf{W}}_i$, $\tilde{\mathbf{U}}_i$ and the $\mathcal{I}$-valued mapping $\tilde{\mathbf{Y}}_i$ by setting

$$\tilde{\mathbf{X}}_i(\mathbf{a}, \mathbf{s}_1, \ldots, \mathbf{s}_n) := \sum_{j \in \mathcal{H}_i \cup \mathcal{L}_i} \left( \frac{\mathbf{a}_j}{\sigma_{i,j}} - \frac{\mathbf{s}_{i,j}}{\sigma_{i,j}} \circ \left( \frac{\alpha_j \mathbf{id}}{\sigma_{i,j}} \right) \right),$$

$$\tilde{\mathbf{Y}}_i(\mathbf{a}, \mathbf{s}_1, \ldots, \mathbf{s}_n) := \sup\left( -\tilde{\mathbf{X}}_i(\mathbf{a}, \mathbf{s}_1, \ldots, \mathbf{s}_n) + \sum_{h=1}^{i-1} G_{i,h} \tilde{\mathbf{W}}_i(\mathbf{a}, \mathbf{s}_1, \ldots, \mathbf{s}_n) \right),$$

$$\tilde{\mathbf{W}}_i(\mathbf{a}, \mathbf{s}_1, \ldots, \mathbf{s}_n) := \tilde{\mathbf{X}}_i(\mathbf{a}, \mathbf{s}_1, \ldots, \mathbf{s}_n)$$
$$- \sum_{h=1}^{i-1} G_{i,h} \tilde{\mathbf{W}}_h(\mathbf{a}, \mathbf{s}_1, \ldots, \mathbf{s}_n) + \tilde{\mathbf{Y}}_i(\mathbf{a}, \mathbf{s}_1, \ldots, \mathbf{s}_n),$$

$$\tilde{\mathbf{U}}_i(\mathbf{a}, \mathbf{s}_1, \ldots, \mathbf{s}_n) := \sum_{j \in \mathcal{H}_i} \left( \frac{\mathbf{s}_{i,j}}{\sigma_{i,j}} \circ \left( \frac{\alpha_j \mathbf{id}}{\sigma_{i,j}} \right) - \frac{\mathbf{a}_j}{\sigma_{i,j}} + \sum_{h=1}^{i-1} \frac{\alpha^{\mathcal{L}}_{h,j} \tilde{\mathbf{W}}_h(\mathbf{a}, \mathbf{s}_1, \ldots, \mathbf{s}_n)}{\sigma_{i,j}} \right).$$

These mappings are well defined for argument function vectors $\mathbf{a} \in \mathcal{D}^{\mathcal{M}}$, $\mathbf{s}_1 \in \mathcal{D}^{\mathcal{H}_1 \cup \mathcal{L}_1}$, ..., $\mathbf{s}_n \in \mathcal{D}^{\mathcal{H}_n \cup \mathcal{L}_n}$ satisfying

$$(21) \quad \sum_{j \in \mathcal{H}_i \cup \mathcal{L}_i} \left( \frac{\underline{\mathbf{s}}_{i,j}}{\sigma_{i,j}} \frac{\alpha_j}{\sigma_{i,j}} - \frac{\overline{\mathbf{a}}_j}{\sigma_{i,j}} \right) > 0 \qquad \text{for every } i \leq n.$$

THEOREM 2. *We let $(\hat{\mathbf{a}}_k)_{k \in \mathbb{N}}$ be a sequence of cumulative arrival processes in $\mathcal{I}^{\mathcal{M}}$ and $(\hat{\mathbf{s}}_{k,i})_{k \in \mathbb{N}}$ a sequence of cumulative service time processes in $\mathcal{I}^{\mathcal{H}_i \cup \mathcal{L}_i}$ for every $i \leq n$. We assume that, for every $k \in \mathbb{N}$, there are vectors $\hat{\alpha}_k \in \mathbb{R}_+^{\mathcal{M}}$, $\hat{\sigma}_{k,1} \in \mathbb{R}_+^{\mathcal{H}_1 \cup \mathcal{L}_1}$, ..., $\hat{\sigma}_{k,n} \in \mathbb{R}_+^{\mathcal{H}_n \cup \mathcal{L}_n}$ such that, almost surely,*

$$\hat{\alpha}_k = \underline{\hat{\mathbf{a}}}_k = \overline{\hat{\mathbf{a}}}_k,$$
$$\hat{\sigma}_{k,i} = \underline{\hat{\mathbf{s}}}_{k,i} = \overline{\hat{\mathbf{s}}}_{k,i} > 0 \qquad \text{for } i = 1, \ldots, n,$$
$$(22) \quad \hat{\rho}_{k,i} := \sum_{j \in \mathcal{H}_i \cup \mathcal{L}_i} \frac{\hat{\alpha}_{k,j}}{\hat{\sigma}_{k,i,j}} < 1 \qquad \text{for } i = 1, \ldots, n.$$

*We assume that the sequence of centered regular network primitives tuples*

$$(\hat{\mathbf{a}}_k - \hat{\alpha}_k \mathbf{id}, \hat{\mathbf{s}}_{k,1} - \hat{\sigma}_{k,1}\mathbf{id}, \ldots, \hat{\mathbf{s}}_{k,n} - \hat{\sigma}_{k,n}\mathbf{id})_{k \in \mathbb{N}}$$

*satisfies a large deviation principle with normalizing sequence $(b_k)_{k \in \mathbb{N}}$ and good rate function $I$ on $\mathcal{D}_0^{\mathcal{M}} \times \mathcal{D}_0^{\mathcal{H}_1 \cup \mathcal{L}_1} \times \cdots \times \mathcal{D}_0^{\mathcal{H}_n \cup \mathcal{L}_n}$. We assume that this rate function takes the value $\infty$ whenever a component of its argument functions is not continuous.*

*If there is a sequence $(d_k)_{k \in \mathbb{N}}$ in $\mathbb{R}_+$ and vectors $\tilde{\alpha} \in \mathbb{R}^{\mathcal{M}}$, $\tilde{\sigma}_1 \in \mathbb{R}^{\mathcal{H}_1 \cup \mathcal{L}_1}$, $\ldots$, $\tilde{\sigma}_n \in \mathbb{R}^{\mathcal{H}_n \cup \mathcal{L}_n}$ with*

$$\infty = \lim_{k \to \infty} d_k,$$
$$\tilde{\alpha} = \lim_{k \to \infty} (\hat{\alpha}_k - d_k \alpha),$$
$$\tilde{\sigma}_i = \lim_{k \to \infty} (\hat{\sigma}_{k,i} - d_k \sigma_i) \qquad \text{for } i = 1, \ldots, n,$$

(23) $$\tilde{\rho}_i := \sum_{j \in \mathcal{H}_i \cup \mathcal{L}_i} \left( \frac{\tilde{\sigma}_{i,j}}{\sigma_{i,j}} \frac{\alpha_j}{\sigma_{i,j}} - \frac{\tilde{\alpha}_j}{\sigma_{i,j}} \right) > 0 \qquad \text{for } i = 1, \ldots, n,$$

*the following 8 statements are valid with normalizing sequence $(b_k)_{k \in \mathbb{N}}$:*

1. *For every $i \leq n$, the sequence $(d_k \mathbf{U}_i(\hat{\mathbf{a}}_k, \hat{\mathbf{s}}_{k,1}, \ldots, \hat{\mathbf{s}}_{k,n}) - d_k(1 - \hat{\rho}_{k,i}^{\mathcal{H}})\mathbf{id})_{k \in \mathbb{N}}$, where*

$$\hat{\rho}_{k,i}^{\mathcal{H}} := \sum_{j \in \mathcal{H}_i} \frac{\hat{\alpha}_{k,j}}{\hat{\sigma}_{k,i,j}},$$

*satisfies a large deviation principle on $\mathcal{D}_0$ with good rate function*

$$\tilde{\mathbf{u}} \mapsto \inf_{\substack{\tilde{\mathbf{a}} \in \mathcal{D}_0^{\mathcal{M}}, \tilde{\mathbf{s}}_1 \in \mathcal{D}_0^{\mathcal{H}_1 \cup \mathcal{L}_1}, \ldots, \tilde{\mathbf{s}}_n \in \mathcal{D}_0^{\mathcal{H}_n \cup \mathcal{L}_n}, \\ \tilde{\mathbf{u}} = \tilde{\mathbf{U}}_i(\tilde{\mathbf{a}} + \tilde{\alpha}\mathbf{id}, \tilde{\mathbf{s}}_1 + \tilde{\sigma}_1\mathbf{id}, \ldots, \tilde{\mathbf{s}}_n + \tilde{\sigma}_n\mathbf{id}) - \tilde{\rho}_i^{\mathcal{H}}\mathbf{id}}} I(\tilde{\mathbf{a}}, \tilde{\mathbf{s}}_1, \ldots, \tilde{\mathbf{s}}_n),$$

*where*

$$\tilde{\rho}_i^{\mathcal{H}} := \sum_{j \in \mathcal{H}_i} \left( \frac{\tilde{\sigma}_{i,j}}{\sigma_{i,j}} \frac{\alpha_j}{\sigma_{i,j}} - \frac{\tilde{\alpha}_j}{\sigma_{i,j}} \right).$$

2. *For every $i \leq n$, the sequence $(d_k \mathbf{V}_i(\hat{\mathbf{a}}_k, \hat{\mathbf{s}}_{k,1}, \ldots, \hat{\mathbf{s}}_{k,n}))_{k \in \mathbb{N}}$ satisfies a large deviation principle on $\mathcal{D}_0$ with good rate function*

$$\tilde{\mathbf{v}} \mapsto \begin{cases} 0, & \text{if } \tilde{\mathbf{v}} = 0, \\ \infty, & \text{otherwise.} \end{cases}$$

3. *For every $i \leq n$, the sequence $(d_k \mathbf{Y}_i(\hat{\mathbf{a}}_k, \hat{\mathbf{s}}_{k,1}, \ldots, \hat{\mathbf{s}}_{k,n}) - d_k(1 - \hat{\rho}_k)\mathbf{id})_{k \in \mathbb{N}}$ satisfies a large deviation principle on $\mathcal{D}_0$ with good rate function*

$$\tilde{\mathbf{y}} \mapsto \inf_{\substack{\tilde{\mathbf{a}} \in \mathcal{D}_0^{\mathcal{M}}, \tilde{\mathbf{s}}_1 \in \mathcal{D}_0^{\mathcal{H}_1 \cup \mathcal{L}_1}, \ldots, \tilde{\mathbf{s}}_n \in \mathcal{D}_0^{\mathcal{H}_n \cup \mathcal{L}_n}, \\ \tilde{\mathbf{y}} = \tilde{\mathbf{Y}}_i(\tilde{\mathbf{a}} + \tilde{\alpha}\mathbf{id}, \tilde{\mathbf{s}}_1 + \tilde{\sigma}_1\mathbf{id}, \ldots, \tilde{\mathbf{s}}_n + \tilde{\sigma}_n\mathbf{id}) - \tilde{\rho}_i\mathbf{id}}} I(\tilde{\mathbf{a}}, \tilde{\mathbf{s}}_1, \ldots, \tilde{\mathbf{s}}_n).$$

4. For every $i \leq n$, the sequence $(d_k \mathbf{W}_i(\hat{\mathbf{a}}_k, \hat{\mathbf{s}}_{k,1}, \ldots, \hat{\mathbf{s}}_{k,n}))_{k \in \mathbb{N}}$ satisfies a large deviation principle on $\mathcal{D}_0$ with good rate function
$$\tilde{\mathbf{w}} \mapsto \inf_{\substack{\tilde{\mathbf{a}} \in \mathcal{D}_0^{\mathcal{M}}, \tilde{\mathbf{s}}_1 \in \mathcal{D}_0^{\mathcal{H}_1 \cup \mathcal{L}_1}, \ldots, \tilde{\mathbf{s}}_n \in \mathcal{D}_0^{\mathcal{H}_n \cup \mathcal{L}_n}, \\ \tilde{\mathbf{w}} = \tilde{\mathbf{W}}_i(\tilde{\mathbf{a}} + \tilde{\alpha}\mathbf{id}, \tilde{\mathbf{s}}_1 + \tilde{\sigma}_1\mathbf{id}, \ldots, \tilde{\mathbf{s}}_n + \tilde{\sigma}_n\mathbf{id})}} I(\tilde{\mathbf{a}}, \tilde{\mathbf{s}}_1, \ldots, \tilde{\mathbf{s}}_n).$$

5. For every $i \leq n$, the sequence $(\mathbf{D}_i(\hat{\mathbf{a}}_k, \hat{\mathbf{s}}_{k,1}, \ldots, \hat{\mathbf{s}}_{k,n}) - \hat{\alpha}_k \mathbf{id})_{k \in \mathbb{N}}$ satisfies a large deviation principle on $\mathcal{D}_0^{\mathcal{M}}$ with good rate function
$$\tilde{\mathbf{d}} \mapsto \inf_{\substack{\tilde{\mathbf{a}} \in \mathcal{D}_0^{\mathcal{M}}, \tilde{\mathbf{s}}_1 \in \mathcal{D}_0^{\mathcal{H}_1 \cup \mathcal{L}_1}, \ldots, \tilde{\mathbf{s}}_n \in \mathcal{D}_0^{\mathcal{H}_n \cup \mathcal{L}_n}, \\ \tilde{\mathbf{d}} = \tilde{\mathbf{a}} - \alpha_i^{\mathcal{L}} \tilde{\mathbf{W}}_i(\tilde{\mathbf{a}} + \tilde{\alpha}\mathbf{id}, \tilde{\mathbf{s}}_1 + \tilde{\sigma}_1\mathbf{id}, \ldots, \tilde{\mathbf{s}}_n + \tilde{\sigma}_n\mathbf{id})}} I(\tilde{\mathbf{a}}, \tilde{\mathbf{s}}_1, \ldots, \tilde{\mathbf{s}}_n).$$

6. For every $i \leq n$, the sequence $(\mathbf{Q}_i(\hat{\mathbf{a}}_k, \hat{\mathbf{s}}_{k,1}, \ldots, \hat{\mathbf{s}}_{k,n}))_{k \in \mathbb{N}}$ satisfies a large deviation principle on $\mathcal{D}_0^{\mathcal{M}}$ with good rate function
$$\tilde{\mathbf{q}} \mapsto \inf_{\substack{\tilde{\mathbf{a}} \in \mathcal{D}_0^{\mathcal{M}}, \tilde{\mathbf{s}}_1 \in \mathcal{D}_0^{\mathcal{H}_1 \cup \mathcal{L}_1}, \ldots, \tilde{\mathbf{s}}_n \in \mathcal{D}_0^{\mathcal{H}_n \cup \mathcal{L}_n}, \\ \tilde{\mathbf{q}} = \alpha_i^{\mathcal{L}} \tilde{\mathbf{W}}_i(\tilde{\mathbf{a}} + \tilde{\alpha}\mathbf{id}, \tilde{\mathbf{s}}_1 + \tilde{\sigma}_1\mathbf{id}, \ldots, \tilde{\mathbf{s}}_n + \tilde{\sigma}_n\mathbf{id})}} I(\tilde{\mathbf{a}}, \tilde{\mathbf{s}}_1, \ldots, \tilde{\mathbf{s}}_n).$$

7. For every $i \leq n$, the sequence $(d_k \mathbf{Z}_i(\hat{\mathbf{a}}_k, \hat{\mathbf{s}}_{k,1}, \ldots, \hat{\mathbf{s}}_{k,n}) - d_k \mathbf{id})_{k \in \mathbb{N}}$ satisfies a large deviation principle on $\mathcal{D}_0^{\mathcal{M}}$ with good rate function
$$\tilde{\mathbf{z}} \mapsto \inf_{\substack{\tilde{\mathbf{a}} \in \mathcal{D}_0^{\mathcal{M}}, \tilde{\mathbf{s}}_1 \in \mathcal{D}_0^{\mathcal{H}_1 \cup \mathcal{L}_1}, \ldots, \tilde{\mathbf{s}}_n \in \mathcal{D}_0^{\mathcal{H}_n \cup \mathcal{L}_n}, \\ \tilde{\mathbf{z}} = e^{\mathcal{L}} \tilde{\mathbf{W}}_i(\tilde{\mathbf{a}} + \tilde{\alpha}\mathbf{id}, \tilde{\mathbf{s}}_1 + \tilde{\sigma}_1\mathbf{id}, \ldots, \tilde{\mathbf{s}}_n + \tilde{\sigma}_n\mathbf{id})}} I(\tilde{\mathbf{a}}, \tilde{\mathbf{s}}_1, \ldots, \tilde{\mathbf{s}}_n).$$

8. The previous statements are also valid in combination: For instance, for every $i, h \leq n$, the sequence
$$(d_k \mathbf{W}_i(\hat{\mathbf{a}}_k, \hat{\mathbf{s}}_{k,1}, \ldots, \hat{\mathbf{s}}_{k,n}), \mathbf{Q}_h(\hat{\mathbf{a}}_k, \hat{\mathbf{s}}_{k,1}, \ldots, \hat{\mathbf{s}}_{k,n}))_{k \in \mathbb{N}}$$
satisfies a large deviation principle on $\mathcal{D}_0 \times \mathcal{D}_0^{\mathcal{M}}$ with good rate function
$$(\tilde{\mathbf{w}}, \tilde{\mathbf{q}}) \mapsto \inf_{\substack{\tilde{\mathbf{a}} \in \mathcal{D}_0^{\mathcal{M}}, \tilde{\mathbf{s}}_1 \in \mathcal{D}_0^{\mathcal{H}_1 \cup \mathcal{L}_1}, \ldots, \tilde{\mathbf{s}}_n \in \mathcal{D}_0^{\mathcal{H}_n \cup \mathcal{L}_n}, \\ \tilde{\mathbf{w}} = \tilde{\mathbf{W}}_i(\tilde{\mathbf{a}} + \tilde{\alpha}\mathbf{id}, \tilde{\mathbf{s}}_1 + \tilde{\sigma}_1\mathbf{id}, \ldots, \tilde{\mathbf{s}}_n + \tilde{\sigma}_n\mathbf{id}), \\ \tilde{\mathbf{q}} = \alpha_h^{\mathcal{L}} \tilde{\mathbf{W}}_h(\tilde{\mathbf{a}} + \tilde{\alpha}\mathbf{id}, \tilde{\mathbf{s}}_1 + \tilde{\sigma}_1\mathbf{id}, \ldots, \tilde{\mathbf{s}}_n + \tilde{\sigma}_n\mathbf{id})}} I(\tilde{\mathbf{a}}, \tilde{\mathbf{s}}_1, \ldots, \tilde{\mathbf{s}}_n).$$

The statements of the theorem remain valid when the elements of all underlying function spaces are restricted to the time interval $\mathbb{R}_+$ and to those starting in 0, and the definitions of the model are modified accordingly. [Hence, conditions (16) and (21) become meaningless.] With this modification, the statements are also true when conditions (22) and (23) are dropped, and/or the topology is changed to the topology of uniform convergence on compacts.

The remarks on Theorem 1 similarly apply to each queue individually. We make the following additional remarks concerning this theorem:

REMARK 5. We say that a pair $(\mathbf{w}, \mathbf{y}) \in \mathcal{D}^n \times \mathcal{I}^n$ *solves the Skorokhod problem* for the vector $\mathbf{z} \in \mathcal{D}^n$ and the matrix $R \in \mathbb{R}^{n \times n}$ (on the time interval $\mathbb{R}$) when the following conditions are met for every component $i \leq n$: $\mathbf{w}_i \geq 0$, $\underline{\mathbf{w}}_i = 0$, $\mathbf{w}_i = \mathbf{z}_i + \sum_{h=1}^n R_{i,h} \mathbf{y}_h$, and $\int_{\mathbb{R}} \mathbf{w}_i(t)\, d\mathbf{y}_i(t) = 0$; compare [9, 14, 20].

Here we define the matrix $R$ as the inverse of the matrix $I + G$, where $I$ is the $n \times n$-identity matrix. Due to the feedforward structure of the queueing networks considered in this work, the matrix $R$ is triangular and the matrix consisting of the absolute values of the entries of the matrix $I - R$ has spectral radius 0. Hence, Theorem 6 in [20] states that, for every $\mathbf{z}$ in the set

$$\mathcal{R} := \{\mathbf{z} \in \mathcal{D}^n : (I+G)\underline{\mathbf{z}} < 0\}$$

(with the usual definition of a matrix-vector product), the Skorokhod problem for $\mathbf{z}$ and $R$ has a unique solution which depends continuously on $\mathbf{z} \in \mathcal{R}$. We denote the map from $\mathbf{z} \in \mathcal{R}$ to the vector of functions $\mathbf{w} \in \mathcal{D}^n$ in the unique solution $(\mathbf{w}, \mathbf{y})$ of the Skorokhod problem for $\mathbf{z}$ and $R$ with $\Phi$. Such a solution map is called a *Skorokhod* or *reflection map*.

By Lemma 1 of [21], we have $\tilde{\mathbf{W}}(\mathbf{a}, \mathbf{s}_1, \ldots, \mathbf{s}_n) = \Phi(R\tilde{\mathbf{X}}(\mathbf{a}, \mathbf{s}_1, \ldots, \mathbf{s}_n))$ for all $\mathbf{a} \in \mathcal{D}^{\mathcal{M}}$, $\mathbf{s}_1 \in \mathcal{D}^{\mathcal{H}_1 \cup \mathcal{L}_1}$, ..., $\mathbf{s}_n \in \mathcal{D}^{\mathcal{H}_n \cup \mathcal{L}_n}$ satisfying condition (22). Under the assumptions of the theorem, we therefore obtain that the sequence of $n$-dimensional processes $(d_k \mathbf{W}(\hat{\mathbf{a}}_k, \hat{\mathbf{s}}_{k,1}, \ldots, \hat{\mathbf{s}}_{k,n}))_{k \in \mathbb{N}}$ satisfies a sample path large deviation principle in $\mathcal{D}_0^n$ with good rate function

$$\tilde{\mathbf{w}} \mapsto \inf_{\substack{\tilde{\mathbf{a}} \in \mathcal{D}_0^{\mathcal{M}}, \tilde{\mathbf{s}}_1 \in \mathcal{D}_0^{\mathcal{H}_1 \cup \mathcal{L}_1}, \ldots, \tilde{\mathbf{s}}_n \in \mathcal{D}_0^{\mathcal{H}_n \cup \mathcal{L}_n}, \\ \tilde{\mathbf{w}} = \Phi(R\tilde{\mathbf{X}}(\tilde{\mathbf{a}} + \tilde{\alpha}\mathbf{id}, \tilde{\mathbf{s}}_1 + \tilde{\sigma}_1\mathbf{id}, \ldots, \tilde{\mathbf{s}}_n + \tilde{\sigma}_n\mathbf{id}))}} I(\tilde{\mathbf{a}}, \tilde{\mathbf{s}}_1, \ldots, \tilde{\mathbf{s}}_n).$$

In this light the rate function for the sequence of $n$-dimensional workload processes in the theorem can be looked upon as a "reflected" rate function in the sense that it is derived from the original rate function through a reflection map. Similarly, one can represent the rate function for queue length and sojourn time processes as "reflected" rate functions.

REMARK 6. When the rate function $I$ in this theorem stems from a moderate deviation principle, it is usually also the rate function in a large deviation principle for the tail probabilities of a Gaussian process. In view of the previous remark and of the fact that the reflection map is positively homogeneous [i.e., $\Phi(c\mathbf{z}) = c\Phi(\mathbf{z})$ for every $c > 0$ and $\mathbf{z} \in \mathcal{R}$], the rate function for the sequence of $n$-dimensional workload processes in this case is also the rate function in a large deviation principle for the tail probabilities of an $n$-dimensional reflected Gaussian process. Similar statements hold for the queue length and sojourn time processes. In this sense our work can rigorously justify the analysis of large deviations of reflected Gaussian processes [1, 20, 23, 25, 27, 28] as a means to analyzing logarithmic tail asymptotics of the behavior of queueing systems in critical loading. Compare Section 5 for an example with renewal and [22] for an example with long-range dependent input processes.

REMARK 7. This theorem parallels the convergence of the queueing processes in distribution to reflected (fractional) Brownian motion in heavy traffic, when the sequence of random network primitives tuples satisfies a functional (fractional) central limit theorem [21, 29].

REMARK 8. Statement 7 of Theorem 2 reveals the validity of the "snapshot principle" [34] in large deviation asymptotics. Namely, the rate function for the sojourn time behavior is only finite, when at every time point the sojourn time of an observer arriving at this time is a specific linear combination of the workload of the nodes at this time.

PROOF OF THEOREM 2. This is analogous to the proof of Theorem 1, but using Theorem 5 of [21] instead of Theorem 4 of [21]. Due to the feedforward structure of the network, one can also prove this theorem by iteratively applying Theorem 1. □

**5. Example: renewal input processes.** In this section we show how Theorem 2 can be applied in a case where the inter arrival and service times of the customers of each class form independent and identically distributed sequences of random variables. We partially parallel Theorem 4.1(c) of [32], but since we aim for the stationary case, and the theory of sample path large deviations in the topology $\|\cdot\|$ is less developed, we need stronger conditions. We show that the obtained rate function is in fact the rate function in a large deviation principle for scaled reflected Brownian motion. When attention is restricted to the positive time interval and initially empty queues, one can use analogous conditions as in Theorem 4.1(c) of [32] and mimic the development of this section in order to obtain large deviation principles for the queueing behavior in the topology of uniform convergence on compacts.

Given a sequence $(X_h)_{h \in \mathbb{Z}}$ of strictly positive, identically distributed and independent random variables with finite mean, we define the linearly interpolated partial sums (lips) process $\mathbf{X}^{\text{lips}} = (\mathbf{X}^{\text{lips}}(t))_{t \in \mathbb{R}}$ by setting, for $h \in \mathbb{Z}$,

$$
(24) \qquad \mathbf{X}^{\text{lips}}(h) := \begin{cases} 0, & \text{if } h = 0, \\ \sum_{m=1}^{h} X_m, & \text{if } h > 0, \\ -\sum_{m=1}^{-h} X_{1-m}, & \text{if } h < 0, \end{cases}
$$

and linearly interpolating these values between subsequent integer time points. Through the law of large numbers we get $\mathbf{X}^{\text{lips}} \in \mathcal{C} \cap \mathcal{I}$ and $\underline{\mathbf{X}^{\text{lips}}} = \overline{\mathbf{X}^{\text{lips}}} = E(X_1)$, with probability one. By modifying the underlying probability space on a zero-set, we can and will assume without loss of generality that these properties are satisfied for all sample paths of linearly interpolated partial sums processes in this work.

Clearly, the time shifted process $\mathbf{X}^{\text{lips}}(\cdot + h) - \mathbf{X}^{\text{lips}}(h)$ has the same distribution as $\mathbf{X}^{\text{lips}}$ for every $h \in \mathbb{Z}$. Furthermore, the process $(\mathbf{X}^{\text{lips}})^{-\mathbf{1}}$ is a renewal process having a jump at time 0 and being linearly interpolated between subsequent jump times.

For $x \in \mathbb{R}$, we let $\lfloor x \rfloor$ be the largest integer $h \in \mathbb{Z}$ with $h \leq x$. If $N$ is an uniformly distributed random variable on $[0, 1[$ independent of $(X_h)_{h \in \mathbb{Z}}$, the process $\hat{\mathbf{x}}$ on $\mathcal{I}$ defined by $\hat{\mathbf{x}}(t) := \lfloor (\mathbf{X}^{\text{lips}} - NX_1)^{-\mathbf{1}}(t) \rfloor$ for $t \in \mathbb{R}$ is a renewal process with stationary increments [i.e., the distribution of $\hat{\mathbf{x}}(\cdot + t) - \hat{\mathbf{x}}(t)$ is the same as the one of $\hat{\mathbf{x}}$ for every $t \in \mathbb{R}$]; see Exercise 1.2.2 in [3].

In the following we shall use these observations to define arrival and service time processes from sequences of random variables. We scale time of the $k$th input process by the factor $k$ and its state by $1/\sqrt{b_k k}$, where $(b_k)_{k \in \mathbb{N}}$ is a sequence of positive numbers satisfying

$$(25) \qquad \lim_{k \to \infty} \frac{b_k}{k} = 0$$

and

$$(26) \qquad \lim_{k \to \infty} \frac{b_k}{\log k} = \infty.$$

We first specify the sequence of renewal arrival processes: For every $k \in \mathbb{N}$ and $j \in \mathcal{M}$, we let $(A_{k,j,h})_{h \in \mathbb{Z}}$ be a sequence of identically distributed strictly positive random variables having finite second moments. For every $k \in \mathbb{N}$ and $j \in \mathcal{M}$, we let $N_{k,j}$ be a random variable which is uniformly distributed on $[0, 1[$. We assume that for every fixed $k \in \mathbb{N}$ the random variables in the set $\{A_{k,j,h} : j \in \mathcal{M}, h \in \mathbb{Z}\} \cup \{N_{k,j} : j \in \mathcal{M}\}$ are independent. We define the process $\hat{\mathbf{a}}_k$ on $\mathcal{I}^{\mathcal{M}}$ by setting, for $j \in \mathcal{M}$ and $t \in \mathbb{R}$,

$$\hat{\mathbf{a}}_{k,j}(t) := \frac{1}{\sqrt{b_k k}} \lfloor (\mathbf{A}_{k,j}^{\text{lips}} - N_{k,j} A_{k,j,1})^{-\mathbf{1}}(kt) \rfloor.$$

Hence, the process $\hat{\mathbf{a}}_{k,j}$ specifies discrete arrivals of class $j$ customers with size $1/\sqrt{b_k k}$ having independent inter arrival times distributed as $A_{k,j,1}/k$. We define, for $k \in \mathbb{N}$ and $j \in \mathcal{M}$,

$$\hat{\alpha}_{k,j} := \sqrt{\frac{k}{b_k}} \frac{1}{E(A_{k,j,1})}.$$

For every $k \in \mathbb{N}$, the process $\hat{\mathbf{a}}_k$ has stationary increments and satisfies

$$\hat{\alpha}_k = \underline{\hat{\mathbf{a}}}_k = \overline{\hat{\mathbf{a}}}_k.$$

Second, we specify service time processes: For every $k \in \mathbb{N}$, $i \leq n$ and $j \in \mathcal{H}_i \cup \mathcal{L}_i$, we let $(S_{k,i,j,h})_{h \in \mathbb{Z}}$ be a sequence of identically distributed strictly positive random variables having finite second moments. We assume that, for every fixed $k \in \mathbb{N}$, the random variables in the set $\{S_{k,i,j,h} : i \leq n, j \in \mathcal{H}_i \cup \mathcal{L}_i, h \in \mathbb{Z}\} \cup \{A_{k,j,h} : j \in \mathcal{M}, h \in \mathbb{Z}\} \cup \{N_{k,j} : j \in \mathcal{M}\}$ are independent. We define the process $\hat{\mathbf{s}}_{k,i}$ on $\mathcal{I}_{\text{inv}}^{\mathcal{H}_i \cup \mathcal{L}_i}$ by setting, for $j \in \mathcal{H}_i \cup \mathcal{L}_i$ and $t \in \mathbb{R}$,

$$\hat{\mathbf{s}}_{k,i,j}(t) := \frac{1}{\sqrt{b_k k}} (\mathbf{S}_{k,i,j}^{\text{lips}})^{-1}(kt).$$

The service time process $\hat{\mathbf{s}}_{k,i,j}$ specifies independent service times distributed as the random variable $S_{k,i,j,1}/k$ for discrete class $j \in \mathcal{H}_i \cup \mathcal{L}_i$ customers of size $1/\sqrt{b_k k}$ at queue $i$. We note that the linear interpolation involved in the definition of $\hat{\mathbf{s}}_{k,i,j}$ does not disrupt this interpretation, since in the definitions of the behavior of the queueing network the process $\hat{\mathbf{s}}_{k,i,j}$ does not appear outside terms of the form $\hat{\mathbf{s}}_{k,i,j}^{-1} \circ \hat{\mathbf{a}}_{k,j}$.

We define, for $k \in \mathbb{N}$, $i \leq n$ and $j \in \mathcal{H}_i \cup \mathcal{L}_i$,

$$\hat{\sigma}_{k,i,j} := \sqrt{\frac{k}{b_k}} \frac{1}{E(S_{k,i,j,1})}.$$

For every $k \in \mathbb{N}$, and $i \leq n$, the process $\hat{\mathbf{s}}_{k,i}$ satisfies

$$\hat{\sigma}_{k,i} = \underline{\hat{\mathbf{s}}}_{k,i} = \overline{\hat{\mathbf{s}}}_{k,i}.$$

The load generated by the sequence of input processes must approach a critical load as $k$ increases. This is captured by assuming the conditions

(27) $$\lim_{k \to \infty} E(A_{k,j,1}) = \alpha_j^{-1} \qquad \text{for } j \in \mathcal{M},$$

(28) $$\lim_{k \to \infty} E(S_{k,i,j,1}) = \sigma_{i,j}^{-1} \qquad \text{for } i \leq n \text{ and } j \in \mathcal{H}_i \cup \mathcal{L}_i,$$

where $\alpha \in \mathbb{R}_+^{\mathcal{M}}$, $\sigma_1 \in \mathbb{R}_+^{\mathcal{H}_1 \cup \mathcal{L}_1}$, ..., $\sigma_n \in \mathbb{R}_+^{\mathcal{H}_n \cup \mathcal{L}_n}$ are "critical rates" vectors with strictly positive components and satisfy conditions (17) and (18). In order to reach the critical load at a specific speed, we additionally assume that there are vectors $\tilde{\alpha} \in \mathbb{R}^{\mathcal{M}}$, $\tilde{\sigma}_1 \in \mathbb{R}^{\mathcal{H}_1 \cup \mathcal{L}_1}$, ..., $\tilde{\sigma}_n \in \mathbb{R}^{\mathcal{H}_n \cup \mathcal{L}_n}$ satisfying

(29) $$\tilde{\alpha}_j = \lim_{k \to \infty} \sqrt{\frac{k}{b_k}} \left( \frac{1}{E(A_{k,j,1})} - \alpha_j \right) \qquad \text{for } j \in \mathcal{M},$$

(30) $$\tilde{\sigma}_{i,j} = \lim_{k \to \infty} \sqrt{\frac{k}{b_k}} \left( \frac{1}{E(S_{k,j,i,1})} - \sigma_{i,j} \right) \qquad \text{for } i \leq n \text{ and } j \in \mathcal{H}_i \cup \mathcal{L}_i.$$

Next, we define a rate function for the sequence of input processes. We assume that there exist vectors $u \in \mathbb{R}_+^{\mathcal{M}}$, $v_1 \in \mathbb{R}_+^{\mathcal{H}_1 \cup \mathcal{L}_1}$, ..., $v_n \in \mathbb{R}_+^{\mathcal{H}_n \cup \mathcal{L}_n}$ with

$$\lim_{k \to \infty} \text{Var}(A_{k,j,1}) = u_j^2 \qquad \text{for } j \in \mathcal{M},$$

$$\lim_{k \to \infty} \text{Var}(S_{k,i,j,1}) = v_{i,j}^2 \qquad \text{for } i \leq n \text{ and } j \in \mathcal{H}_i \cup \mathcal{L}_i.$$

We define the good rate function $I^{\text{Brown}} \colon \mathcal{D}_0 \to \mathbb{R}_+ \cup \{\infty\}$ for $\mathbf{x} \in \mathcal{D}_0$ by

$$I^{\text{Brown}}(\mathbf{x}) := \begin{cases} \int_{\mathbb{R}} \dfrac{\dot{\mathbf{x}}(t)^2}{2}\, dt, & \text{if } \mathbf{x}(0) = 0, \text{ and } \mathbf{x} \text{ is absolutely continuous,} \\ \infty, & \text{otherwise.} \end{cases}$$

Here $\dot{\mathbf{x}}$ denotes a derivative of an absolutely continuous function $\mathbf{x} \in \mathcal{D}_0$. This rate function is well known from Schilder's theorem [8] dealing with logarithmic tail asymptotics of Brownian motion. We let $I^{\text{renewal}}$ be the good rate function on $\mathcal{D}_0^{\mathcal{M}} \times \mathcal{D}_0^{\mathcal{H}_1 \cup \mathcal{M}_1} \times \cdots \times \mathcal{D}_0^{\mathcal{H}_n \cup \mathcal{M}_n}$ given by

$$I^{\text{renewal}}(\tilde{\mathbf{a}}, \tilde{\mathbf{s}}_1, \ldots, \tilde{\mathbf{s}}_n) := \sum_{j \in \mathcal{M}} \frac{\alpha_j^3 I^{\text{Brown}}(\tilde{\mathbf{a}}_j)}{u_j^2} + \sum_{i=1}^n \sum_{j \in \mathcal{H}_i \cup \mathcal{L}_i} \frac{\sigma_{i,j}^3 I^{\text{Brown}}(\tilde{\mathbf{s}}_{i,j})}{v_{i,j}^2}.$$

Here we set $x/0 := 0$, if $x = 0$, and $x/0 := \infty$, if $x \in \mathbb{R} \setminus \{0\}$, in order to capture the case where some of the limiting variances are zero. We can now apply Corollaries A.2 and A.3 of Appendix A in order to get the following:

COROLLARY 3. *If there exists $\delta > 0$ and $c > 0$ such that, for every $k \in \mathbb{N}$ and $y \in [-\delta, \delta]$,*

$$\log E(\exp(y(A_{k,j,1} - E(A_{k,j,1})))) \leq cy^2 \qquad \text{for } j \in \mathcal{M},$$

$$\log E(\exp(y(S_{k,i,j,1} - E(S_{k,i,j,1})))) \leq cy^2 \qquad \text{for } i \leq n \text{ and } j \in \mathcal{H}_i \cup \mathcal{L}_i,$$

*the sequence*

$$(\hat{\mathbf{a}}_k - \hat{\alpha}_k \mathbf{id}, \hat{\mathbf{s}}_{k,1} - \hat{\sigma}_{k,1} \mathbf{id}, \ldots, \hat{\mathbf{s}}_{k,n} - \hat{\sigma}_{k,n} \mathbf{id})_{k \in \mathbb{N}}$$

*satisfies a sample path large deviation principle on $\mathcal{D}_0^{\mathcal{M}} \times \mathcal{D}_0^{\mathcal{H}_1 \cup \mathcal{L}_1} \times \cdots \times \mathcal{D}_0^{\mathcal{H}_n \cup \mathcal{L}_n}$ with good rate function $I^{\text{renewal}}$ and normalizing sequence $(b_k)_{k \in \mathbb{N}}$.*

If for a $k \in \mathbb{N}$ the stability condition (22) is satisfied, the random network primitives tuple $(\hat{\mathbf{a}}_k, \hat{\mathbf{s}}_{k,1}, \ldots, \hat{\mathbf{s}}_{k,n})$ is regular, and the queueing processes $\mathbf{U}_i(\hat{\mathbf{a}}_k, \hat{\mathbf{s}}_{k,1}, \ldots, \hat{\mathbf{s}}_{k,n})$, $\mathbf{V}_i(\hat{\mathbf{a}}_k, \hat{\mathbf{s}}_{k,1}, \ldots, \hat{\mathbf{s}}_{k,n})$, $\mathbf{Y}_i(\hat{\mathbf{a}}_k, \hat{\mathbf{s}}_{k,1}, \ldots, \hat{\mathbf{s}}_{k,n})$, $\mathbf{W}_i(\hat{\mathbf{a}}_k, \hat{\mathbf{s}}_{k,1}, \ldots, \hat{\mathbf{s}}_{k,n})$, $\mathbf{D}_i(\hat{\mathbf{a}}_k, \hat{\mathbf{s}}_{k,1}, \ldots, \hat{\mathbf{s}}_{k,n})$, $\mathbf{Q}_i(\hat{\mathbf{a}}_k, \hat{\mathbf{s}}_{k,1}, \ldots, \hat{\mathbf{s}}_{k,n})$ and $\mathbf{Z}_i(\hat{\mathbf{a}}_k, \hat{\mathbf{s}}_{k,1}, \ldots, \hat{\mathbf{s}}_{k,n})$ are well defined for every queue $i \leq n$. From Lemma B.2 of Appendix B, we deduce that the processes $\mathbf{V}(\hat{\mathbf{a}}_k, \hat{\mathbf{s}}_{k,1}, \ldots, \hat{\mathbf{s}}_{k,n})$, $\mathbf{W}(\hat{\mathbf{a}}_k, \hat{\mathbf{s}}_{k,1}, \ldots, \hat{\mathbf{s}}_{k,n})$, $\mathbf{Q}(\hat{\mathbf{a}}_k, \hat{\mathbf{s}}_{k,1}, \ldots, \hat{\mathbf{s}}_{k,n})$ and $\mathbf{Z}(\hat{\mathbf{a}}_k, \hat{\mathbf{s}}_{k,1}, \ldots, \hat{\mathbf{s}}_{k,n}) - \mathbf{id}$ are jointly stationary.

COROLLARY 4. *If, for every $k \in \mathbb{N}$, conditions (22) and (23), and the assumptions of Corollary 3 are satisfied, the prerequisites and consequences of Theorem 2 are valid with $I := I^{\text{renewal}}$, and $d_k := \sqrt{k/b_k}$ for $k \in \mathbb{N}$.*

Clearly, independent exponentially or deterministically distributed interarrival and service times with mean values satisfying conditions (22), (23), (27), (28), (29) and (30) constitute a possible choice through which the assumptions of the previous two corollaries are satisfied.

For a positive semidefinite covariance matrix $V \in \mathbb{R}^{n \times n}$ and a vector $b \in R^n$, we define the good rate function $I_V : \mathcal{D}_0^n \to \mathbb{R} \cup \{\infty\}$ by

$$I_V^{\text{Brown}}(\mathbf{x}) := \begin{cases} \int_{\mathbb{R}} \frac{\dot{\mathbf{x}}(t)^T V^{-1} \dot{\mathbf{x}}(t)}{2} \, dt, & \text{if } \mathbf{x} \text{ is absolutely continuous} \\ & \text{and } \mathbf{x}(0) = 0, \\ \infty, & \text{otherwise.} \end{cases}$$

Here $\dot{\mathbf{x}}(t)$ is the vector of derivatives of the absolutely continuous function $\mathbf{x} \in \mathcal{D}^n$ at time $t \in \mathbb{R}$, and $\dot{\mathbf{x}}(t)^T$ its transpose. If the matrix $V$ is degenerate, the expression $\dot{\mathbf{x}}(t)^T V^{-1} \dot{\mathbf{x}}(t)/2$ in this definition must be understood as $\sup_{y \in \mathbb{R}^n}(y^T \dot{\mathbf{x}}(t) - y^T V y / 2)$; compare Remark 4.1 in [32]. We note that if $\hat{\mathbf{B}}$ is an $n$-dimensional Brownian motion with stationary increments on the time interval $\mathbb{R}$, passing through 0 at time 0, having covariance matrix $V$ at time 1, and possessing drift 0, the sequence $(\hat{\mathbf{B}}/\sqrt{b_k})_{k \in \mathbb{N}}$ satisfies a large deviation principle with normalizing sequence $(b_k)_{k \in \mathbb{N}}$ good rate function $I_V^{\text{Brown}}$ on $\mathcal{D}_0^n$ [8].

Here we define $V := RUR^T$, where $R^T$ is the transpose of the matrix $R$ and

$$U_{i,h} := \begin{cases} \displaystyle\sum_{j \in (\mathcal{H}_i \cup \mathcal{L}_i) \cap (\mathcal{H}_h \cup \mathcal{L}_h)} \frac{u_j^2}{\alpha_j^3 \sigma_{i,j} \sigma_{h,j}}, & \text{if } i \neq h, \\ \displaystyle\sum_{j \in \mathcal{H}_i \cup \mathcal{L}_i} \left( \frac{u_j^2}{\alpha_j^3 \sigma_{i,j}^2} + \frac{v_{i,j}^2 \alpha_j}{\sigma_{i,j}^6} \right), & \text{if } i = h. \end{cases}$$

Furthermore, we set $\tilde{\zeta} := R\tilde{\rho}$, where the vector $\tilde{\rho} > 0$ has been defined in (22). Then Remark 5 after Theorem 2 together with the contraction principle [7] and simple algebraic calculations yields that, under the assumptions of Corollary 4, the sequence $(d_k \mathbf{W}(\hat{\mathbf{a}}_k, \hat{\mathbf{s}}_{k,1}, \ldots, \hat{\mathbf{s}}_{k,n}))_{k \in \mathbb{N}}$ satisfies a sample path large deviation principle on $\mathcal{D}_0^n$ with normalizing sequence $(b_k)_{k \in \mathbb{N}}$ and good rate function

$$\tilde{\mathbf{w}} \mapsto \inf_{\tilde{\mathbf{z}} \in \mathcal{D}_0^n, \, \tilde{\mathbf{w}} = \Phi(\tilde{\mathbf{z}} - \tilde{\zeta}\mathbf{id})} I_V^{\text{Brown}}(\tilde{\mathbf{z}}).$$

The contraction principle shows that this rate function is also the rate function in the large deviation principle for the sequence of stationary reflected Brownian motions $(\Phi(\hat{\mathbf{B}}/\sqrt{b_k} - \tilde{\zeta}\mathbf{id}))_{k \in \mathbb{N}}$ with normalizing sequence $(b_k)_{k \in \mathbb{N}}$. In particular, results concerning the rate function of reflected Brownian motion [2, 10, 15, 16, 18, 24] might be applied to obtain statements about the asymptotic behavior of tail probabilities of multiclass feedforward queueing networks in heavy traffic.

## APPENDIX A: MODERATE DEVIATIONS FOR RENEWAL PROCESSES

Here we consider for every $k \in \mathbb{N}$ a sequence $(X_{k,h})_{h \in \mathbb{Z}}$ of independent and identically distributed strictly positive random variables with finite second moments. We derive moderate deviation principles for certain sequences of renewal processes constructed from these random variables in the topology $\|\cdot\|$, from known moderate deviation principles for partial sums processes in the topology of uniform convergence on compacts ([33], Lemma 6.1).

LEMMA A.1. *If there are values $\xi > 0$ and $\Sigma \geq 0$ such that*
$$\lim_{k \to \infty} E(X_{k,1}) = \xi,$$
$$\lim_{k \to \infty} \mathrm{Var}(X_{k,1}) = \Sigma^2,$$
*and constants $c > 0$ and $\delta > 0$ such that, for every $y \in [-\delta, \delta]$ and $k \in \mathbb{N}$,*
$$\Lambda_k(y) := \log E(y(X_{k,1} - E(X_{k,1}))) \leq cy^2,$$
*then for every sequence $(b_k)_{k \in \mathbb{N}}$ of positive numbers satisfying (25) and (26), the centered and scaled sequence of linearly interpolated partial sums processes*
$$\left( \frac{\mathbf{X}_k^{\mathrm{lips}}(k \cdot) - E(X_{k,1})\mathbf{id}(k \cdot)}{\sqrt{b_k k}} \right)_{k \in \mathbb{N}}$$
*satisfies a sample path large deviation principle on $\mathcal{D}_0$ with good rate function $I^{\mathrm{Brown}}/\Sigma^2$ and normalizing sequence $(b_k)_{k \in \mathbb{N}}$.*

PROOF. We define for $k \in \mathbb{N}$ the piecewise constant partial sums (pcps) process $\mathbf{X}_k^{\mathrm{pcps}}$ by setting $\mathbf{X}_k^{\mathrm{pcps}}(t) := \mathbf{X}_k^{\mathrm{lips}}(\lfloor t \rfloor)$ for $t \in \mathbb{R}$. Under the assumptions of the lemma, Example 7.2 of [30] yields that the sequence of processes
$$\left( \frac{\mathbf{X}_k^{\mathrm{pcps}}(k \cdot) - E(X_{k,1})\mathbf{id}(k \cdot)}{\sqrt{b_k k}} \right)_{k \in \mathbb{N}}$$
satisfies a sample path large deviation principle on $\mathcal{D}_0$ with good rate function $I^{\mathrm{Brown}}/\Sigma^2$ and normalizing sequence $(b_k)_{k \in \mathbb{N}}$ in the topology of uniform convergence on compacts.

For every $\varepsilon > 0$ and $\tau > 0$, we have, for all $k \in \mathbb{N}$ with $\varepsilon \sqrt{b_k k} \geq 2 \sup_{k \in \mathbb{N}} E(X_{k,1})$,
$$P\left( \sup_{-k\tau - 1 \leq h \leq k\tau + 1} \frac{X_{k,h}}{\sqrt{b_k k}} > \varepsilon \right)$$
$$\leq (2k\tau + 2) P(X_{k,1} > \varepsilon \sqrt{b_k k})$$
$$\leq (2k\tau + 2) P(\exp(\delta(X_{k,1} - E(X_{k,1})) - \delta\varepsilon \sqrt{b_k k}))$$
$$\leq (2k\tau + 2) \exp(c\delta^2 - \delta\varepsilon \sqrt{b_k k}).$$

Hence, for every $\varepsilon > 0$,

$$\limsup_{k \to \infty} \frac{1}{b_k} \log P\left(\sup_{t \in [-\tau, \tau]} \frac{|\mathbf{X}_k^{\text{lips}}(kt) - \mathbf{X}_k^{\text{pcps}}(kt)|}{\sqrt{b_k k}} > \varepsilon\right)$$

$$\leq \limsup_{k \to \infty} \frac{1}{b_k} \log P\left(\sup_{-k\tau - 1 \leq h \leq k\tau + 1} \frac{X_{k,h}}{\sqrt{b_k k}} > \varepsilon\right)$$

$$\leq \limsup_{k \to \infty} \left(\frac{\log(2k\tau + 2) + c\delta^2}{b_k} - \delta\sqrt{\frac{k}{b_k}}\right) = -\infty.$$

This implies by Theorem 4.2.13 in [7] that the sequence of processes

$$\left(\frac{\mathbf{X}_k^{\text{lips}}(k \cdot) - E(X_{k,1})\mathbf{id}(k \cdot)}{\sqrt{b_k k}}\right)_{k \in \mathbb{N}}$$

also satisfies a sample path large deviation principle on $\mathcal{D}_0$ with good rate function $I^{\text{Brown}}/\Sigma^2$ and normalizing sequence $(b_k)_{k \in \mathbb{N}}$ in the topology of uniform convergence on compacts.

Next, we strengthen this large deviation principle to the topology induced by the norm $\|\cdot\|$. According to Theorem 18 in [20], one can reach this goal by showing that, for every $\varepsilon, \nu > 0$, there exists a $\tau \geq 0$ such that

$$\limsup_{k \to \infty} \frac{1}{b_k} \log P\left(\sup_{|t| \geq \tau} \frac{|\mathbf{X}_k^{\text{lips}}(kt) - E(X_{k,1})kt|}{t\sqrt{b_k k}} > \varepsilon\right) < -\nu.$$

We choose $\varepsilon, \nu > 0$. A Chernoff type estimate gives, for every $h \in \mathbb{N}$ and $\gamma \in [0, 2c\delta]$,

$$P\left(\sum_{m=1}^{h} (X_{k,m} - E(X_{k,1})) > h\gamma\right)$$

$$\leq E \exp\left(\frac{\gamma}{2c}\left(\sum_{m=1}^{h} (X_{k,m} - E(X_{k,1})) - h\gamma\right)\right)$$

$$\leq E \exp\left(\frac{\gamma}{2c} \sum_{m=1}^{h} (X_{k,m} - E(X_{k,1})) - \frac{h\gamma^2}{2c}\right)$$

$$= \exp\left(h\Lambda_k\left(\frac{\gamma}{2c}\right) - \frac{h\gamma^2}{2c}\right) \leq \exp\left(-\frac{h\gamma^2}{4c}\right).$$

Similarly, one obtains, for every $h \in \mathbb{N}$ and $\gamma \in [0, 2c\delta]$,

$$\log P\left(\sum_{m=1}^{h} (X_{k,m} - E(X_{k,1})) < -h\gamma\right) \leq -\frac{h\gamma^2}{4c}.$$

Thus, we get, for every $\tau \in \mathbb{N}$ and $k \in \mathbb{N}$ with $\varepsilon\sqrt{b_k/k} \leq 2c\delta$,

$$P\left(\sup_{t \geq \tau} \frac{|\mathbf{X}_k^{\text{lips}}(kt) - E(X_{k,1})kt|}{t\sqrt{b_k k}} > \varepsilon\right)$$

$$\leq \sum_{h=k\tau}^{\infty} P\left(\left|\sum_{m=1}^{h}(X_{k,m} - E(X_{k,1}))\right| > h\varepsilon\sqrt{\frac{b_k}{k}}\right)$$

$$\leq \sum_{h=k\tau}^{\infty} 2\exp\left(-\frac{hb_k\varepsilon^2}{4ck}\right) = \frac{2\exp(-\tau b_k\varepsilon^2/(4c))}{1-\exp(-b_k\varepsilon^2/(4ck))}.$$

This implies that, for $\tau \in \mathbb{N}$ with $\tau \geq 4c\nu/\varepsilon^2$,

$$\limsup_{k\to\infty} \frac{1}{b_k} \log P\left(\sup_{t \geq \tau} \frac{|\mathbf{X}_k^{\text{lips}}(kt) - E(X_{k,1})kt|}{t\sqrt{b_k k}} > \varepsilon\right)$$

$$\leq \limsup_{k\to\infty}\left(-\frac{\tau\varepsilon^2}{4c} - \frac{\log(1-\exp(-\varepsilon^2 b_k/(4ck)))}{b_k}\right)$$

$$= -\nu + \limsup_{k\to\infty}\left(\frac{\log(4ck/(\varepsilon^2 b_k))}{b_k}\right) = -\nu.$$

Here we have used the assumptions (25) and (26). A similar bound holds on the negative time interval which yields the envisaged strengthening through the principle of the largest term. □

We define the process $\hat{\mathbf{x}}_k$ by $\hat{\mathbf{x}}_k(t) := (\mathbf{X}_k^{\text{lips}})^{-1}(kt)$ for $t \in \mathbb{R}$. Analogously to Lemma 2.4 in [32], we get the following:

COROLLARY A.2. *Under the assumptions of Lemma A.1, the sequence of processes*

$$\left(\frac{\hat{\mathbf{x}}_k - E(X_{k,1})^{-1}k\mathbf{id}}{\sqrt{b_k k}}\right)_{k \in \mathbb{N}}$$

*satisfies a sample path large deviation principle on $\mathcal{D}_0$ with good rate function $\xi^3 I^{\text{Brown}}/\Sigma^2$ and normalizing sequence $(b_k)_{k\in\mathbb{N}}$.*

For $k \in \mathbb{N}$, we let $N_k$ be a random variable which is independent of the sequence $(X_{k,h})_{h\in\mathbb{Z}}$ and uniformly distributed on $[0,1]$, and define the process $\hat{\mathbf{x}}'_k(t) := \lfloor(\mathbf{X}_k^{\text{lips}} - NX_{k,1})^{-1}\rfloor(kt)$. Clearly, we have $|\hat{\mathbf{x}}'_k - \hat{\mathbf{x}}_k| \leq 1$, which yields by exponential equivalence [7]

COROLLARY A.3. *Under the assumptions of Lemma A.1, the sequence of processes*

$$\left(\frac{\hat{\mathbf{x}}'_k - E(X_{k,1})^{-1}k\mathbf{id}}{\sqrt{b_k k}}\right)_{k \in \mathbb{N}}$$

satisfies a sample path large deviation principle on $\mathcal{D}_0$ with good rate function $\xi^3 I^{\mathrm{Brown}}/\Sigma^2$ and normalizing sequence $(b_k)_{k\in\mathbb{N}}$.

## APPENDIX B: STATIONARITY OF NETWORK BEHAVIOR

In this section we state independence and stationary increments conditions for random regular network primitives which imply the stationarity of the workload and queue length processes.

We first investigate the behavior of a single queueing node under time shifts. For $c \in \mathbb{R}$ and a function $\mathbf{d} \in \mathcal{D}$, we let $\Theta_c \mathbf{d} \in \mathcal{D}$ be the *time shifted function* defined by $(\Theta_c \mathbf{d})(t) := \mathbf{d}(t+c)$ for every $t \in \mathbb{R}$. Furthermore, we set $\Xi_c \mathbf{d} := \Theta_c \mathbf{d} - \mathbf{d}(c)$, which implies $(\Xi_c \mathbf{d})(0) = 0$ for every $c \in \mathbb{R}$ and $\mathbf{d} \in \mathcal{D}$. We use these shift mappings on product function spaces by applying them componentwise. For instance, if $c \in \mathbb{R}^n$, $\mathbf{d} \in \mathcal{D}^n$ and $i \leq n$, we let the $i$th component of $\Xi_c \mathbf{d} \in \mathcal{D}^n$ be defined by $(\Xi_c \mathbf{d})_i := \Xi_{c_i} \mathbf{d}_i$.

LEMMA B.1. *If $(\mathbf{a}, \mathbf{s}) \in \mathcal{I}^\mathcal{M} \times \mathcal{I}_{\mathrm{inv}}^{\mathcal{H} \cup \mathcal{L}}$ is a regular node primitives pair and $t \in \mathbb{R}$,*

$$\Xi_{\mathbf{a}(t)}(\mathbf{s}^{-1}) = (\Xi_{\mathbf{s}^{-1}\circ\mathbf{a}(t)}\mathbf{s})^{-1},$$

$$\Xi_t(\mathbf{s}^{-1}\circ\mathbf{a}) = (\Xi_{\mathbf{s}^{-1}\circ\mathbf{a}(t)}\mathbf{s})^{-1}\circ(\Xi_t\mathbf{a}),$$

$$\Theta_t \mathbf{V}^\mathcal{H}(\mathbf{a},\mathbf{s}) = \mathbf{V}^\mathcal{H}(\Xi_t\mathbf{a}, \Xi_{\mathbf{s}^{-1}\circ\mathbf{a}(t)}\mathbf{s}),$$

$$\Theta_t \mathbf{W}^{\mathcal{H},\mathcal{L}}(\mathbf{a},\mathbf{s}) = \mathbf{W}^{\mathcal{H},\mathcal{L}}(\Xi_t\mathbf{a}, \Xi_{\mathbf{s}^{-1}\circ\mathbf{a}(t)}\mathbf{s}),$$

$$\Theta_t \mathbf{Q}^{\mathcal{H},\mathcal{L}}(\mathbf{a},\mathbf{s}) = \mathbf{Q}^{\mathcal{H},\mathcal{L}}(\Xi_t\mathbf{a}, \Xi_{\mathbf{s}^{-1}\circ\mathbf{a}(t)}\mathbf{s}),$$

$$\Theta_t(\mathbf{Z}^{\mathcal{H},\mathcal{L}}(\mathbf{a},\mathbf{s}) - \mathbf{id}) = \mathbf{Z}^{\mathcal{H},\mathcal{L}}(\Xi_t\mathbf{a}, \Xi_{\mathbf{s}^{-1}\circ\mathbf{a}(t)}\mathbf{s}) - \mathbf{id}$$

*and*

$$\Xi_t \mathbf{D}^{\mathcal{H},\mathcal{L}}(\mathbf{a},\mathbf{s}) = \mathbf{D}^{\mathcal{H},\mathcal{L}}(\Xi_t\mathbf{a}, \Xi_{\mathbf{s}^{-1}\circ\mathbf{a}(t)}\mathbf{s}) + \mathbf{Q}^{\mathcal{H},\mathcal{L}}(\Xi_t\mathbf{a}, \Xi_{\mathbf{s}^{-1}\circ\mathbf{a}(t)}\mathbf{s})(0).$$

PROOF. Under the assumptions of the lemma, the first two lines in the display follow from the fact that, for every $b, c \in \mathbb{R}$ and $\mathbf{d} \in \mathcal{I}_{\mathrm{inv}}$, we have $(\Theta_c \mathbf{d} - b)^{-1} = \Theta_b(\mathbf{d}^{-1}) - c$. Next, one obtains the third and fourth line of the display by using the equation $\Theta_c(\sup \mathbf{d}) - b = \sup(\Theta_c \mathbf{d} - b)$, which is valid for every $\mathbf{d} \in \mathcal{D}_{\sup}$ and $c, b \in \mathbb{R}$. For $j \in \mathcal{H}$ and $c, t \in \mathbb{R}$, we calculate

$$-(\Theta_c \mathbf{Q}_j^{\mathcal{H},\mathcal{L}}(\mathbf{a},\mathbf{s}))(t) = \mathbf{D}_j^{\mathcal{H},\mathcal{L}}(\mathbf{a},\mathbf{s})(t+c) - \mathbf{a}_j(t+c)$$

$$= \sup_{\tau \leq t+c,\ \sum_{\ell \in \mathcal{H}} \mathbf{s}_\ell^{-1}\circ\mathbf{a}_\ell(\tau) \leq \sum_{\ell \in \mathcal{H}} \mathbf{s}_\ell^{-1}\circ\mathbf{a}_\ell(t+c) - \mathbf{V}^\mathcal{H}(\mathbf{a},\mathbf{s})(t+c)} \mathbf{a}_j(\tau) - \mathbf{a}_j(t+c)$$

$$= \sup_{\tau - c \leq t,\ \sum_{\ell \in \mathcal{H}}(\Xi_c(\mathbf{s}_l^{-1}\circ\mathbf{a}_\ell))(\tau-c) \leq \sum_{\ell \in \mathcal{H}}(\Xi_c(\mathbf{s}_\ell^{-1}\circ\mathbf{a}_\ell))(t) - \Theta_c \mathbf{V}^\mathcal{H}(\mathbf{a},\mathbf{s})(t)} (\Xi_c \mathbf{a}_j)(\tau - c)$$

$$- (\Xi_c \mathbf{a}_j)(t)$$
$$= \mathbf{D}_j^{\mathcal{H},\mathcal{L}}(\Xi_c \mathbf{a}, \Xi_{\mathbf{s}^{-1} \circ \mathbf{a}(c)} \mathbf{s})(t) - (\Xi_c \mathbf{a}_j)(t)$$
$$= -\mathbf{Q}_j^{\mathcal{H},\mathcal{L}}(\Xi_c \mathbf{a}, \Xi_{\mathbf{s}^{-1} \circ \mathbf{a}(c)} \mathbf{s})(t).$$

This implies the fifth equation for components $j \in \mathcal{H}$. The remaining equations follow along similar transformations which are omitted. $\square$

We use these shift properties to identify situations where the behavior of the feedforward network is stationary.

LEMMA B.2. *We let* $(\hat{\mathbf{a}}, \hat{\mathbf{s}}_1, \ldots, \hat{\mathbf{s}}_n) \in \mathcal{I}^{\mathcal{M}} \times \mathcal{I}_{\mathrm{inv}}^{\mathcal{H}_1 \cup \mathcal{L}_1} \times \cdots \times \mathcal{I}_{\mathrm{inv}}^{\mathcal{H}_n \cup \mathcal{L}_n}$ *be a random regular network primitives tuple satisfying the following four conditions:*

1. *The process* $\hat{\mathbf{a}}$ *is independent of* $\hat{\mathbf{s}}_i$ *for every* $i \leq n$.
2. *The processes* $\hat{\mathbf{s}}_{i,j}$ *for* $i \leq n$ *and* $j \in \mathcal{H}_i \cup \mathcal{L}_i$ *are independent.*
3. *The process* $\hat{\mathbf{a}}$ *has stationary increments, that is, the process* $\Xi_t \hat{\mathbf{a}}$ *has the same distribution for every* $t \in \mathbb{R}$.
4. *For every* $i \leq n$ *and* $j \in \mathcal{H}_i \cup \mathcal{L}_i$, *the process* $\hat{\mathbf{s}}_{i,j}^{-1}$ *has stationary increments for time shifts with values from the image of* $\hat{\mathbf{a}}_j$, *which means that the distribution of the process* $\Xi_a(\hat{\mathbf{s}}_{i,j}^{-1})$ *is the same for every* $a \in \{\hat{\mathbf{a}}_j(t) \colon t \in \mathbb{R}\} \subset \mathbb{R}$.

*Then the processes* $\mathbf{V}_1(\hat{\mathbf{a}}, \hat{\mathbf{s}}_1, \ldots, \hat{\mathbf{s}}_n), \ldots, \mathbf{V}_n(\hat{\mathbf{a}}, \hat{\mathbf{s}}_1, \ldots, \hat{\mathbf{s}}_n), \mathbf{W}_1(\hat{\mathbf{a}}, \hat{\mathbf{s}}_1, \ldots, \hat{\mathbf{s}}_n),$ $\ldots, \mathbf{W}_n(\hat{\mathbf{a}}, \hat{\mathbf{s}}_1, \ldots, \hat{\mathbf{s}}_n), \mathbf{Q}_1(\hat{\mathbf{a}}, \hat{\mathbf{s}}_1, \ldots, \hat{\mathbf{s}}_n), \ldots, \mathbf{Q}_n(\hat{\mathbf{a}}, \hat{\mathbf{s}}_1, \ldots, \hat{\mathbf{s}}_n),$ *and* $\mathbf{Z}_1(\hat{\mathbf{a}}, \hat{\mathbf{s}}_1, \ldots, \hat{\mathbf{s}}_n) - \mathbf{id}, \ldots, \mathbf{Z}_n(\hat{\mathbf{a}}, \hat{\mathbf{s}}_1, \ldots, \hat{\mathbf{s}}_n) - \mathbf{id}$ *are jointly stationary, that is, their common distribution is invariant under the shift* $\Theta_t$ *for every* $t \in \mathbb{R}$.

PROOF. The conditions of the lemma together with the second equation of Lemma B.1 imply that the distribution of the pair $(\Xi_t \hat{\mathbf{a}}, \Xi_{\hat{\mathbf{s}}_1^{-1} \circ \hat{\mathbf{a}}(t)} \hat{\mathbf{s}}_1)$ does not depend on $t \in \mathbb{R}$. Hence, the joint stationarity of the processes $\mathbf{V}_1(\hat{\mathbf{a}}, \hat{\mathbf{s}}_1, \ldots, \hat{\mathbf{s}}_n), \mathbf{W}_1(\hat{\mathbf{a}}, \hat{\mathbf{s}}_1, \ldots, \hat{\mathbf{s}}_n), \mathbf{Q}_1(\hat{\mathbf{a}}, \hat{\mathbf{s}}_1, \ldots, \hat{\mathbf{s}}_n)$ and $\mathbf{Z}_1(\hat{\mathbf{a}}, \hat{\mathbf{s}}_1, \ldots, \hat{\mathbf{s}}_n) - \mathbf{id}$ follows from the third to sixth equation of lemma B.1 and the definitions of these processes.

Setting $\hat{\mathbf{d}}_1 := \mathbf{D}_1(\hat{\mathbf{a}}, \hat{\mathbf{s}}_1, \ldots, \hat{\mathbf{s}}_n) = \mathbf{D}^{\mathcal{H}_1, \mathcal{L}_1}(\hat{\mathbf{a}}, \hat{\mathbf{s}}_1)$, the last equation of Lemma B.1 yields that the distribution of the triple $(\Xi_t \hat{\mathbf{a}}, \Xi_{\hat{\mathbf{s}}_1^{-1} \circ \hat{\mathbf{a}}_1(t)} \hat{\mathbf{s}}_1, \Xi_t \hat{\mathbf{d}}_1)$ does not depend on $t \in \mathbb{R}$. Clearly, the image of $\hat{\mathbf{d}}_{1,j}$ is a subset of the image of $\hat{\mathbf{a}}_j$ for every $j \in \mathcal{M}$. In view of the independence conditions 1 and 2, the process $\hat{\mathbf{d}}_1$ is independent of the processes $\hat{\mathbf{s}}_2, \ldots, \hat{\mathbf{s}}_n$. In particular, by Assumption 4, the distribution of the quadruple $(\Xi_t \hat{\mathbf{a}}, \Xi_{\hat{\mathbf{s}}_1^{-1} \circ \hat{\mathbf{a}}_1(t)} \hat{\mathbf{s}}_1, \Xi_t \hat{\mathbf{d}}_1, \Xi_{\hat{\mathbf{s}}_2^{-1} \circ \hat{\mathbf{d}}_1(t)} \hat{\mathbf{s}}_2)$ does not depend on $t \in \mathbb{R}$. Repeating the reasoning of the first paragraph simultaneously for the first two queues, one therefore obtains the joint stationarity of the considered queueing processes of queues 1 and 2.

Iterating the scheme of the second paragraph through the remaining queues yields the statement of the lemma. $\square$

SIEMENS AG
CT SE 6
81730 MÜNCHEN
GERMANY
E-MAIL: Kurt.Majewski@siemens.com